\documentclass[12pt]{amsart}
\usepackage{amsmath}
\usepackage{amssymb}
\usepackage{amscd}
\usepackage{graphicx}

\newtheorem*{TA}{\sc theorem A}
\newtheorem*{TAP}{\sc theorem A'}

\newtheorem{thm}{\sc theorem}[section]
\newtheorem{pro}[thm]{\sc proposition}
\newtheorem{lem}[thm]{\sc lemma}
\newtheorem{cor}[thm]{\sc corollary}

\theoremstyle{definition}

\newtheorem{dfn}[thm]{\sc definition}

\theoremstyle{remark}

\newcommand{\Int}{\mbox{$\mbox{\rm Int}$}}

\newcommand{\diff}{\mbox{$\mbox{\rm Diff}$}}

\newcommand{\pr}{\mbox{$\mbox{\rm pr}$}}

\def\R{{\bf R}}
\def\C{{\bf C}}

\def\S{{\bf S}}

\def\D{{\bf D}}
\def\T{{\bf T}}

\def\F{{\mathcal F}}
\def\NN{{\mathcal N}}

\def\H{{\mathcal H}}
\def\SS{{\mathcal S}}

\def\X{{\mathcal X}}

\def\SO{{\rm SO}}

\def\id{{\rm id}}

\def\mun{{^{-1}}}

\def\G{{\Gamma_1}}

\def\half{{1\over 2}}

\def\and{{\rm and}}
\def\min{{\rm min}}
\def\round{{\rm round}}
\def\holes{{\rm holed}}
\def\whirl{{\rm whirl}}

\def\can{{\rm can}}

\def\plug{{\rm plug}}

\def\sing{{\rm Sing}}

\begin{document}

\title{Regularization and minimization of Haefliger structures
of codimension one.}

\author{ Ga\"el Meigniez}
\address{{UMR 6205, Laboratoire de Math\'ematiques de Bretagne Atlantique,
  Universit\'e de Bretagne-Sud, Universit\'e Europ\'eenne de Bretagne ---}
Postal address: Universit\'e de Bretagne-Sud, L.M.B.A.,
BP 573, F-56017 Vannes, France}
\email{Gael.Meigniez@univ-ubs.fr}

\date\today



\dedicatory{Dedicated to Paul Schweitzer.}

\begin{abstract}
In dimensions 4 and more, a proof is given
of Thurs\-ton's h-principle for foliations of codimension one,
by elementary means
not using Mather's homology equivalence. Moreover,
the produced foliations are minimal, that is, all leaves
are dense. This makes minimal smooth foliations
on every closed manifold of dimension at least 4
whose Euler characteristic is null.

\end{abstract}
\maketitle

\section{Introduction}

Thuston's {h-principle} for foliations of codimension one
\cite{thurston76} is a powerful method to
 produce such foliations on a given compact manifold $M^n$~. It allows
to prescribe both the homotopy class of the foliation regarded
as a hyperplane field, and as a $\G$-structure.

See paragraph \ref{CVN} below for the
 Haefliger structures of codimension one,
also called
$\G$-structures. Practically, one may think
of them as the singular foliations of codimension one
 with Morse-type singularities; which is the generic case. A homotopy between
two $\G$-structures on $M$ is a $\G$-structure on $M\times[0,1]$~.

In this introduction, to fix ideas, one reduces oneself to the
case where $M$ is closed; to the
smooth ($C^\infty$) differentiability class; and
 to the co-oriented
case: all hyperplane fields, foliations, and $\G$-structures, are understood co-oriented. See section \ref{generalizations_sec} for the general case.

The h-principle can be stated in terms of a \emph{formal foliation} on $M$~,
by which one means a pair $(\xi,
\tau)$ where $\xi$ is a $\G$-structure, and $\tau$ is a hyperplane
field on $M$~. Two formal foliations $(\xi_0,\tau_0)$~,
$(\xi_1,\tau_1)$ are \emph{homotopic} if
$\xi_0$ is homotopic to $\xi_1$ (as $\G$-structures)
and if, independently,
 $\tau_0$ is homotopic to $\tau_1$ (as hyperplane fields).
A foliation is nothing but a formal foliation that is holonomic:
$\xi$ is regular and $\tau$ is tangent to $\xi$~.
\begin{thm} (Thurston \cite{thurston76})
 On every closed manifold,
every formal foliation
 is 
homotopic to some
foliation.
\end{thm}

 Both proofs given in \cite{thurston76}
 follow the same scheme of construction in three steps.
 First one makes on $M$ a foliation with ``holes'', that is, parts of the manifold
 left unfoliated. Second, one fills the holes, that is,
one extends the foliation into them.
 The filling argument is substantial, although
it is constructive and its means are completely elementary.
Third, to adjust the homotopy class of the foliation
as a Haefliger structure, one changes it, in some domains of $M$~,
 for
other pieces of foliations, given by Mather's homology equivalence \cite{mather1}
between, on the one hand, the classifying space $B\diff_c(\R)^\delta$ of the
group of
diffeomorphisms with compact supports, and on the other hand,
the loop space
 of Haefliger's classifying space $B\G^+$ for the co-oriented
$\G$-structures.

The produced foliation has compact leaves, because while filling
the holes some kind of Reeb components are created inside,
 and because in general
there are also some kind of Reeb components in the
 pieces given by Mather's homology equivalence.
\medbreak

It is notorious \cite{thurston94} that the subsequent
 study of fol\-iations was strained by a dram\-atic
 sociological phenomenon,
induced by the feeling that the
achievements of the homotopic theory
of foliations in the first half of the 70's
 were
difficult; and by the belief that
 Thurston had "cleaned out"
and even "killed" the subject  (as a doctorate student, I was
still told so in Paris in the middle 80's). In
particular, the paradigm grew
 that in codimension one,
 the use of Mather's homology equivalence
 was unavoidable to prove the h-principle, despite the fact that in 1976 Thurston himself did not write so.

 {\it ``\emph{Currently,} step 3 \emph{seems} to involve some nonelementary background about classifying spaces for
 Haefliger structures and groups of diffeomorphisms.''} \cite{thurston76}
\medbreak
The first aim of the present paper is to give, in dimensions $n\ge 4$~, a constructive
proof of this h-principle by elementary means,
not involving such background.

Moreover, we make the h-principle produce foliations which are \emph{
minimal}, that is, all leaves are dense in $M$~.
\medbreak

\begin{TA}
 On every closed manifold of dimension at least $4$~,
every formal foliation
 is 
homotopic to some \emph{minimal}
foliation.
\end{TA}
\begin{cor} Every closed, connected manifold of
dimension at least $4$ and whose Euler characteristic is null,
admits a minimal, smooth foliation of codimension one.
\end{cor}
For example, $\S^5$ does, as well as $\S^3\times\S^2$~, $\S^3\times\S^1$~, and $(\S^2\times\S^2)\sharp
(\S^3\times\S^1)\sharp(\S^3\times\S^1)$~. This corollary answers a classical
 question raised by Lawson, and
 disproves
a classical conjecture about an analogue, in high dimensions,
 of Novikov's closed leaf theorem \cite{lawson}.

\medbreak
Recall that P. Schweitzer proved that on every
closed manifold of dimension at least $4$~,
every formal foliation can be homotoped to some
foliation of class $C^1$ without compact leaf \cite{schweitzer}.

In dimension $3$~,
a constructive proof of Thurston's h-principle for codimension-one
foliations by elementary means, not using Mather's homology
equivalence,
was given recently, prior to the present work \cite{lm}. The
produced foliation has a precise geometric form, linked to some
open book decomposition of the manifold.
\medbreak

\medbreak
Our proof of theorem A works directly on the Morse singular foliation
$\xi$ to be regularized and minimized through homotopy. Its singularities are cancelled by pairs,
in Morse's way. Some holes, isomorphic to Thurston's, appear,
as the price to pay to put the
 pairs of singularities into cancellation position. The holes are filled without interior
compact leaves. This allows minimality. Some care in
the choice of the cancelling pairs and in the cancellation itself,
allows to prescribe the homotopy class of the produced foliation,
regarded as a hyperplane field.

Here is a more precise account. Given a $\G$-structure $\xi$ and a nonsingular
 vector field $V$ on $M$~, a first generic homotopy
turns $\xi$ into a Morse singular foliation. Then,
pairs of Morse singularities of indices 1 and 2 are created
 to give some
genus to the leaves, after what the leaves are made dense by some
local whirls.
This first step presents no difficulty.

After some new homotopies of $\xi$ that
create new singularities by pairs, and after a convenient homotopy
of $V$~, the Morse singular foliation $\xi$
admits a pseudogradient $\nabla\xi$ such that the points of $M$
where $\nabla\xi$ is opposed to $V$ constitute a family of arcs
transverse to $\xi$~, each arc joining two singularities
of $\xi$ of successive indices.

Then, a homotopy of $\xi$ in a neighborhood of each arc places
the pair of singularities into cancellation position, in Morse's
 sense; at the price of a hole. Then, one cancels all pairs of singularities,
leaving on $M$ a foliation with holes.
Thanks to the very choice of the arcs,
and to some
Thom-Pontryagin like considerations,
the gradient of the
 produced foliation falls naturally into
the homotopy class of $V$~, as a nonsingular vector field.

Finally one fills each hole. Here, the minimality
of the surrounding foliation allows to enlarge, following
\cite{thurston76}, the
hole by tunnelling out a \emph{worm gallery} through the foliation, before filling the enlarged hole.
Our holes are exactly the same
as Thurston's, and our filling argument is in the largest part
borrowed from his, but with two critical differences.
We avoid to create kind of Reeb components
in the hole,
and
so, the produced foliation remains minimal. While "rolling up
the holes", we use $\S^{k-1}\times\S^1$ where
\cite{thurston76} uses the $k$-torus. This second difference
simplifies the inductive "rolling up" process,
and allows that
the produced
foliation stay in the homotopy class of $\xi$~, as a $\G$-structure. The filling argument remains substantial.

The paper is organized in a thematic and progressive way.
 In section \ref{minimize_sec}, one produces minimal foliations. In section
\ref{homotopy_sec}, one controls their homotopy classes,
 as $\G$-structures.
 In section
\ref{hyperplane_sec}, one controls their homotopy classes, as hyperplane fields.
Finally, section \ref{generalizations_sec} gives the theorem
 in its full generality.

\medbreak It is a pleasure to thank Fran\c cois Laudenbach for
 sharing his interrogations about the actual status of Mather's
homology equivalence
 in Thurston's construction of codimension one foliations;
 and for his listening and observations all along this work.

\subsection{Conventions, vocabulary and notations.}\label{CVN}

The object under consideration is a co-oriented
 {\it Haefliger structure} of codimension one --- more briefly a \emph{$\G$-structure} \cite{haefliger} --- on $M$~.
 It can be defined as
 a foliation $\xi$
of codimension one on some neighborhood of $M\times 0$
in $M\times\R$~, transverse to each fibre $x\times\R$~;
or more precisely, as  the germ along $M\times 0$
 of such a foliation.

The \emph{singularities} of $\xi$ are the points of $M$ where $M\times 0$ is not transverse to $\xi$~. Write $\sing(\xi)\subset M$
their set. The restriction of $\xi$ to $M\times 0$ induces
 a foliation on $M\setminus\sing(\xi)$~.
 One makes no difference
between ``foliation'' and ``regular $\G$-structure''.

A \emph{homotopy} between two $\G$-structures $\xi_0$~, $\xi_1$ on $M$~, is a $\G$-structure on the manifold $M\times[0,1]$
whose restriction to $M\times i$ is $\xi_i$~, for $i=0,1$~.

 Every real-valued function $f$ \emph{defines} on its domain
 a co-oriented $\G$-structure $\xi$~: the
pullback through $f$ of the regular $\G$-structure
on the real line. One also
 calls $f$ a \emph{first integral} of $\xi$.
Obviously, every $\G$-structure admits
a local first integral in a neighborhood of every point.
\medbreak
One denotes $\D^k$ the compact unit ball in $\R^k$~,
 and $\S^{k-1}:=\partial\D^k$~. In
particular, $\D^1$ is the interval $[-1,+1]$~.
A basepoint is generally written $*$~.
Also, one sometimes regards $\D^2$ as the unit disk in $\C$~,
and uses the basepoint $1$ in $\S^1=\partial \D^2$~.

 In any product $X\times Y$~, where $Y$ is an interval or a circle,
 the
projection to the second factor is called the \emph{height function.}
 The foliation that it defines is the \emph{height foliation.} The unit vector field
 negatively parallel to $Y$ is the \emph{height gradient.}

By \emph{whirling} a foliation $\F$
 in a domain diffeomorphic to some product
 $X\times\D^1$ where $\F$
is the height foliation, we mean changing $\F$ in this domain
for the suspension of some representation
of $\pi_1X$ into the diffeomorphisms of $\D^1$~.

The \emph{stabilization} of a foliated manifold $(V,\F)$ by a manifold $X$
is the foliated manifold $(X\times V,\pr_2^*\F)$~.

\section{Making a minimal foliation.}\label{minimize_sec}

Let $M$ be a closed connected
manifold of dimension $n\ge 4$ whose Euler characteristic is null.
The object of this section is to endow $M$ with
 a minimal foliation of codimension $1$ and class $C^\infty$~.
 To the next sections
is postponed the care of its homotopy classes, as
a $\G$-structure and as a hyperplane field.

 In this section, the smooth (that is, $C^\infty$) differentiability class
 is understood; and all $\G$-structures
 are understood co-oriented.

 The construction will be in three steps. 

A $\G$-structure is qualified \emph{Morse} iff its local first integrals (see \ref{CVN}) are Morse functions.
We also call it for short a \emph{Morse foliation.}
(In the sequel, by the term "foliation" without the qualificative "Morse" we continue to understand a \emph{regular} $\G$-structure.)

First one will make a minimal Morse foliation on $M$~. Then the
singularities will be cancelled by pairs, at the price of
some domains in $M$ left unfoliated --- so-called \emph{holes}.
Third, the holes will be filled, that is,
 the foliation will be extended inside them.

\medbreak
Our tools will just be the creation and the cancellation of
singularities by pairs, following Morse
 \cite{milnor}\cite{morse}. The conditions of
the cancellation need to be reviewed from a viewpoint
that fits Morse foliations. Let $\xi$ be a Morse foliation on $M$~.

On the unit disk
 $\D^k$, the \emph{canonical} $\G$-structure $\xi^k_\can$
 is the one defined by the
square of the norm function.

\begin{dfn}\label{stable_disk} A
 \emph{stable disk} for
 $\xi$ at its singularity $s$ of index $i$~, is
an embedding of $(\D^{n-i},0)$ into $(M,s)$
 through which the pullback of
 $\xi$ equals $\xi_\can^{n-i}$~, and their co-orientations match.

An \emph{unstable disk} for
 $\xi$ at $s$~, is
an embedding of $(\D^{i},0)$ into $(M,s)$
 through which the pullback of
 $\xi$  equals $\xi_\can^{i}$~, and their co-orientations
are opposed.
\end{dfn}
So, at $s$~, the hessian of
$\xi$ is positive (resp. negative) definite on the linear
space tangent to the stable (resp. unstable) disk.

\begin{dfn}\label{cancel_dfn} A \emph{cancellation} pair of disks $D$~, $D'$
 for $\xi$ at its singularities $s$ and $s'$~, of respective indices
$i$ and $i+1$ (where $0\le i\le n-1$)~, is as follows:
\begin{enumerate}
\item $D$ is a stable disk
at $s$ and $D'$ is an unstable disk at $s'$~;
\item Their boundaries meet in a single point $x$~,
at which $\partial D$
and $\partial D'$  are transverse in the leaf of $\xi$
 through $x$~;
\item Otherwise $D$ and $D'$ are disjoint.
\end{enumerate}
\end{dfn}

\begin{lem}\label{cancel_lem} {\rm (Cancellation)}
 Let $M$~,
 $\xi$~, $s$~, $s'$~, $D$~, $D'$ be as above; and let $N$
be a small enough neighborhood of $D\cup D'$ in $M$~.

Then, $N$
 admits a foliation $\NN$
equal to $\xi$ outside some compact subset in $N$~.
Moreover, no leaf of $\NN$
is relatively compact in $N$~.
\end{lem}
\begin{proof}
Obviously, if $N$ is small enough, $\xi$ admits in $N$
 a first integral $f$~, null
 in restriction to the bouquet of spheres $\partial D\cup\partial D'$~.
This bouquet has a compact neighborhood
$E\subset N$ diffeomorphic to $B\times[-\epsilon,
+\epsilon]$~, where $B$ is
 $\S^{n-i-1}\times\S^i$ minus
an open $(n-1)$-disk, and such that
 $f\vert E$ is the height function.

Then, a move classical in Morse theory,
changing $f$ on some smaller neighborhood of $D$~,
makes the critical value at $s$ lift up to $-\epsilon/2$~.
Symmetrically, one makes
 the critical value at $s'$ descend down to $\epsilon/2$
 (the critical values don't cross each other).
One gets on $N$ a Morse function $f'$~, equal to $f$ outside
some compact subset, and whose singularities
are still $s$ and $s'$~. The bouquet has
a compact neighborhood $E'\subset N$ diffeomorphic to
$\D^{n-1}\times[-\epsilon,+\epsilon]$~,
containing $s$ and $s'$ in its interior, and such that
$f'\vert E'$ is the height function on a neighborhood of
 $\partial E'$~.

($E'$ is obtained from $E$ by attaching an
 $(i+1)$-handle along $*\times\S^i\times(+\epsilon)$
and a $(n-i)$-handle along $\S^{n-i-1}\times *\times(-\epsilon)$~, both attachment spheres being
 endowed with the trivial framing.)

One defines $\NN$ by
 the height function on $E'$
and by $f'$ on
$N\setminus E'$~. Obviously, $\NN$ has
no leaf relatively compact in $N$~.
\end{proof}

\subsection{First step: making a minimal Morse foliation.}
\label{minimize_sbs}
The construction of a minimal Morse foliation on every manifold
of dimension at least $3$ is not difficult,
and already known \cite{bonatti}.
Here we give a construction which fits the needs of section
\ref{homotopy_sec}.

One starts with any Morse foliation $\xi$ on $M$~,
 for example the one defined by some Morse function.

Then, one gets
 rid of local extrema, that is,
 singularities of extremal indices $0$ and $n$~.
To this aim, close to each singularity $s$ of $\xi$ of index $0$~,
one creates a pair of singularities $s'$~, $s''$ of
respective
 indices $1$~, $2$~.
Then, $s$ and $s'$ admit respectively
 a stable $n$-disk $D$ and an unstable $1$-disk $D'$
(definition \ref{stable_disk}), such
that their intersection is a single extremity of $D'$~, contained in
 $\partial D$ (notice that the two extremities of a stable or
 unstable $1$-disk
do not lie necessarily in the same leaf) . The
cancellation pair of disks $D$~, $D'$
allows
the cancellation of $s$
and $s'$ (lemma \ref{cancel_lem}).
As a whole, $s$ has been changed for $s''$~, which is of index 2.
We have also created an $n$-dimensional Reeb component.

In the same way, every singularity of index $n$
is changed for one of index $n-2$~, and a Reeb component.
\medbreak

Then, $\xi$ being free of local extrema,
it admits a total transversal, that is,
 a finite disjoint union of compact arcs $A_j$ embedded in $M$ transversely
to $\xi$~, and such that every leaf meets at least one of them.
One will make $\xi$ minimal by the following modification
 in some small
 neighborhood $N_j$ of each arc $A_j$~.
 
Each $N_j$ is identified with
$\D^{n-1}\times\D^1$~, in such a way that $A_j=0\times[-1/4,+1/4]$~; and
that $\xi\vert N_j$ is the height foliation (see \ref{CVN}).
A Morse function $f$ is made
 on
 $N_j$~, from the height function,
 by creating two pairs of singularities $(s_1, s_2)$
and $(s'_1,s'_2)$
 of respective
indices 1, 2, 1, 2, in cancellation position.
 The function $f$ is still the height function close to $\partial N_j$~. One makes
$f(s_1)$~, $f(s'_1)<-1/2$ and $f(s_2)$~, $f(s'_2)>1/2$~.
 The domain $f\mun[-1/2,+1/2]$
is diffeomorphic to $L\times[-1/2,+1/2]$~,
 where $L$ is the $(n-1)$-disk with two 1-handles.

Two diffeomorphisms $\alpha$~, $\alpha'$ of the interval $[-1/2,+1/2]$~,
 both flat on the identity at $\pm1/2$~,
 are chosen such that every orbit
 in the open interval
 $(-1/2,+1/2)$ under the group generated by these diffeomorphisms
 is dense there. One has the representation
 $$\pi_1 L\to \diff_+[-1/2,+1/2]:\gamma\mapsto(\alpha)^{(a\gamma)}(\alpha')^{(a'\gamma)}$$
 where $a$ and $a'$ are Poincar\'e-dual to the stable attachment spheres of $s_1$
 and $s'_1$~, respectively.

Then we whirl $f$ between the singularities.
That is, we define a Morse foliation
 $\xi_\plug$ on $N_j$
by $f$~, except in the domain $f\mun[-1/2,+1/2]$~,
where $\xi_\plug$ is the suspension of the above representation.
 One can arrange that  $\xi_\plug$ coincides
with the height foliation close to
 $\partial N_j$~.
Every leaf of $\xi_\plug$ in $f\mun(-1/2,+1/2)$ is dense there.

Having changed $\xi$ to $\xi_\plug$ inside each $N_j$~, we obtain a new Morse foliation $\xi_\min$ on $M$~.
Every leaf of $\xi_\min$ meets some solid cylinder $N_j$ in its open subset
 $f\mun(-1/2,+1/2)$~, and is dense there.
So, every leaf of $\xi_\min$ is locally dense.
So, every leaf of $\xi_\min$ is dense in $M$~.

\subsection{Second step: canceling
 the singularities.}\label{cancel_sbs}
The second step of the construction of a minimal foliation on $M$
is to cancel the singularities of the minimal Morse foliation
 $\xi_\min$~, at the price of leaving some
 unfoliated
holes.

\medbreak
In order to apply the cancellation lemma \ref{cancel_lem}, one
first arranges that the singularities can be matched into pairs of
successive indices, without asking for the moment any cancellation position.
To get such a partition of $\sing(\xi_\min)$ into pairs,
 one just creates, anywhere in $M$~, some
pairs of new Morse singularities of appropriate successive indices,
as follows.

Recall
that $\xi_\min$ has no singularities of index $0$ or $n$~,
 and we won't create any.
 First,
creating if necessary some pairs of indices $2$
and $3$~, arrange that there are at least as many singularities
of index $2$ as of index $1$~. Then
match every singularity of index $1$ with one of index $2$~. Then,
creating if necessary some pairs of indices $3$
and $4$~, arrange that there are at least as many singularities
of index $3$ as unmatched singularities of index $2$~. Then
match every unmatched singularity of index $2$ with one of index $3$~.
 And so on, until all singularities remaining unmatched are of index
$n-2$ or $n-1$~. Since the Euler characteristic of $M$ is null,
they are in the same number, and the matching of
$\sing(\xi_\min)$ is complete. (This argument,
 and further down
 the existence of transverse arcs, are actually the only global
arguments in this paper.) 

\begin{figure}
\includegraphics*[scale=0.45, angle=-90]{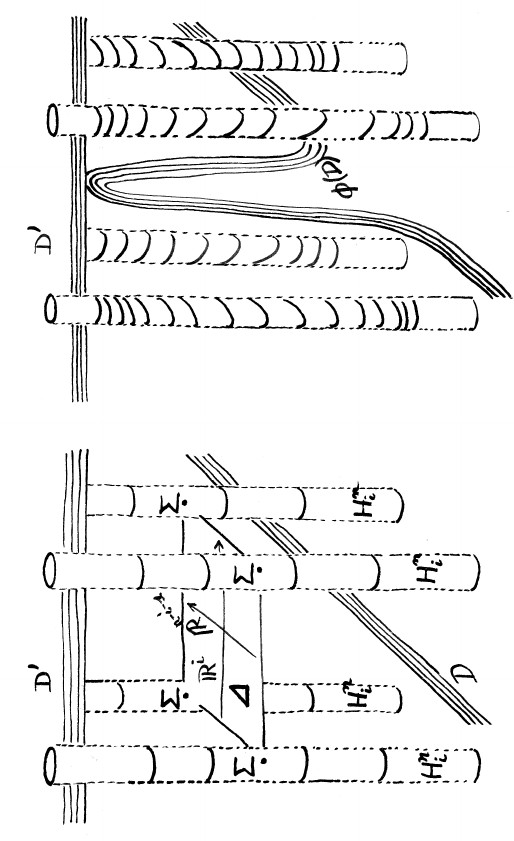}
\caption{Before and after the whirl.}
\label{whirl_fig}
\end{figure}

Still denote $\xi_\min$ the resulting minimal Morse foliation.
\medbreak
Consider one of the pairs $s$~, $s'$ of matched singularities.
Their indices are respectively $i\ge 1$ and $i+1\le n-1$~.
Obviously, at every Morse singularity there exist a small stable
 disk and a small unstable one (definition \ref{stable_disk}).
Let $D$ be a small stable $(n-i)$-disk at $s$ and $D'$ be a small unstable $(i+1)$-disk
at $s'$~, disjoint from $D$~. The leaves of $\xi_\min$ being dense,
 $M$ contains an embedded arc $A$~, positively
transverse to $\xi$~, from some boundary point of $D$ to
some boundary point of $D'$~, and otherwise disjoint from $D$ and $D'$~.
We are going to force the condition of lemma \ref{cancel_lem}, that is,
the existence of a pair of cancellation disks intersecting once on their boundaries, by
a whirl (see paragraph \ref{CVN}) of $\xi_\min$ close to $A$~,
at the price of a hole.
\medbreak
A small open neighborhood $U$ of $A$ is identified with$$\R^n=\R^{n-i-1}\times\R^{i}\times\R$$
in such a way that
$\xi_\min\vert U$ is the height foliation, defined by
 the last coordinate $x_n$~; and that
$$D\cap U=\R^{n-i-1}\times 0\times(-\infty,-1/2]$$
$$D'\cap U=0\times\R^i\times[1/2,+\infty)$$
 (figure \ref{whirl_fig}, left).
One considers in $\R^{n-1}=\R^{n-i-1}\times\R^i$ the unit bidisk
and its corners:
$$\Sigma:=\S^{n-i-2}\times\S^{i-1}\subset\D^{n-i-1}\times\D^i:=\Delta$$
Every loop $\gamma$ in $\R^{n-1}\setminus\Sigma$ has a linking number $\ell(\gamma)$
 with $\Sigma$~, that is,
the algebraic number of intersections of $\gamma$ with
$\D^{n-i-1}\times\S^{i-1}$~.
On the other hand, one chooses a diffeomorphism $\phi$ of $\R$~,
whose support is compact and contained in $(-1,+1)$~,
and such that $\phi(-1/2)=+1/2$
(in the next paragraph,
more conditions will be imposed on $\phi$~, to be able
to fill the holes.)
One gets
a representation
$$\rho:\pi_1(\R^{n-1}\setminus\Sigma)\to
\diff_c(\R):\gamma\mapsto\phi^{\ell(\gamma)}$$
the suspension of which gives a foliation $\xi_\whirl$ on $\R^{n}\setminus(\Sigma\times\D^1)$~,
 well-defined up to a vertical isotopy, ``vertical'' meaning parallel
to the $x_n$-axis.
We can arrange
that $\xi_\whirl$ is the height foliation on $\R^n\setminus(\Delta\times\D^1)$
and on
 $0\times\R^i\times\R$~.
In restriction to $\R^{n-i-1}\times 0\times\R$~,
the foliation $\xi_\whirl$ is the image of
the height foliation by 
some vertical, compactly supported isotopy $\Phi$~.
 The boundaries
$\partial(D'\cap U)$ and $\partial\Phi(D\cap U)$ meet at
 the point $(0,\dots,0,1/2)$~, transversely in the leaf
of $\xi_\whirl$ through this point.

\medbreak
The set that remains unfoliated is $\Sigma\times\D^1$~, but
it is better to regard
as unfoliated the interior of its small tubular neighborhood, called a \emph{hole}:
$$H_i^n\cong\D^2\times\Sigma\times\D^1$$ 
(figure \ref{whirl_fig}, right).

Both components $\D^2\times\Sigma\times\pm 1$ of $\partial H_i^n$
are tangent to $\xi_\whirl$ and without holonomy. On the
contrary, $(\partial\D^2)\times\Sigma\times\D^1$ is transverse
to $\xi_\whirl$~, and the restriction of $\xi_\whirl$
to $(\partial\D^2)\times\Sigma\times\D^1$ endows it
with a structure of foliated $\D^1$-bundle above the base
$(\partial\D^2)\times\Sigma$~, whose
monodromy above $\partial\D^2$ is $\phi$ (due to the whirl),
while the monodromy above any loop in $\Sigma$ is the identity
of $\D^1$~.
\medbreak
At the price of leaving this hole unfoliated, $s$ can be cancelled with $s'$~.
Let us think again of $U\cong\R^n$ as an open subset in $M$~.
Extending $\Phi$ to the all of $D$
 by the identity,
one gets an embedding $\Phi$ of $D$ into $M\setminus H^n_i$~,
isotopic to the identity: the small
 disk $D$ pushes a pseudopod along
the arc $A$ until it touches the boundary of $D'$~.
 Extending $\xi_\whirl$
by $\xi_\min$ outside $U$~,
 one gets a Morse foliation on $M\setminus H^n_i$~, for which $\Phi(D)$ and $D'$ form a cancellation pair
 (definition \ref{cancel_dfn}).
Cancel $s$ with $s'$ (lemma \ref{cancel_lem}).

Having done so for each pair of matched singularities, one gets
a (regular) foliation $\xi_\holes$ on $M_\holes:=M$ minus the holes.

Finally, $\xi_\holes$ is also minimal, in
the sense that all of its leaves are dense in $M_\holes$~. Indeed,
clearly  the minimality property of $\xi_\min$ was preserved through the whirl, and through
the cancellation of the singularities (thanks to the last sentence of lemma \ref{cancel_lem}).

\subsection{Third step: filling the holes.}\label{filling_sbs}
Here much is borrowed from Thurston \cite{thurston76}; but there are also substantial differences.
A complete account is necessary
as well for the seek of sections \ref{homotopy_sec} and
 \ref{hyperplane_sec}.
 I shall indicate the places where we depart significatively from \cite{thurston76}.

In a general way, by a \emph{hole} we mean
a compact manifold $H$
 with boundary
 and
corners, endowed along $\partial H$ with a germ $\H$ of foliation
 of codimension one in $H$~.
\emph{Filling} the hole means extending this germ into a
 foliation $\F$ inside $H$~. We are also interested to fill it \emph{without
interior leaf:} every leaf of $\F$ meets
$\partial H$~.

Most holes will have the following \emph{standard} form. One has  $H=B\times Y$~,
where $Y=\D^1$ or $\S^1$~, and where $B$ is
a compact connected manifold
with nonempty connected smooth boundary.
On some neighborhood of $B\times\partial Y$~, the germ $\H$
is the height foliation.
 On the other hand,
 $\H\vert(\partial B)\times Y$
 is a foliated $Y$-bundle, that is,
the suspension of some representation
$$\rho:\pi_1\partial B\to\diff_c(-1,+1){\ \rm or\ }\diff_+(\S^1)$$
For every loop $\gamma$ in $\partial B$~,
we also refer to $\rho(\gamma)$ as the \emph{monodromy}
of the hole above $\gamma$~.
 We call a standard hole \emph{straight} if
$Y=\D^1$~, \emph{round} if $Y=\S^1$~, \emph{discal} if its base
splits as a product
$B=\D^2\times X$~,
where $X$ is a closed manifold.
Here are the relevant examples of standard holes, and their
diagnostics.
\medbreak
(i) In case $\rho$ extends to some representation
from $\pi_1B$ to $\diff_c(-1,+1)$ or to $\diff_+(\S^1)$~,
 then of course we fill $H$ by the suspension
 of the extended representation, and this filling is without interior leaves.

(ii) By contrast,
every $\phi\in \diff_c(-1,+1)$ defines a 3-dimensional straight
discal hole $\D^2\times\D^1$ whose monodromy above $\partial\D^2$ is $\phi$~.
By Reeb's global stability theorem, this standard hole \emph{cannot}
be filled unless $\phi$ is the identity.

(iii) More generally, consider a
straight discal hole $\D^2\times X\times\D^1$~, where $X$ is a closed manifold;
the monodromy being trivial above $\pi_1X$,
and nontrivial above the loop $(\partial\D^2)\times *$~, which is
compressible in $B$~. In other words, $H$ is the $X$-stabilization (see \ref{CVN}) of the example (ii).
Then, by an easy generalization of Reeb's
global stability theorem,
$H$ cannot be filled (exercise). Unfortunately, every hole $H^n_i$ left by paragraph \ref{cancel_sbs}
falls to this case. Following Thurston, one will enlarge $H_i^n$
by a worm gallery (see further down) and then be able to fill the enlarged hole.

(iv) Actually, the method to fill each enlarged
 $H^n_i$ will be to divide it
into smaller holes that will fall either to the suspension
case (i) above, or to the following one.

 Say that some subset of $\D^1$ \emph{brackets} some other one~,
 if they are disjoint and if
 every point of the second lies between
two points of the first.

Write $\T^r:=(\S^1)^r$ the $r$-torus.

Let $Y:=\D^1$ or $\S^1$~, and $r\ge 1$~. Given $r+1$ commuting diffeomorphisms
 $\phi$~, $\psi_1$~, \dots, $\psi_r\in\diff_c(-1,+1)$ (resp. $\diff_+(\S^1)$)~,
consider the $(r+3)$-dimensional discal hole$$H:=\D^2\times\T^r\times Y$$
whose monodromy is $\phi$ over $\partial\D^2$ and $\psi_k$ over the $k$-th $\S^1$ factor ($1\le k\le r$)~.
\begin{lem}\label{bracketting_inessential}
\emph{(Bracketted compressible monodromy)}
Assume that~:
\begin{enumerate}
\item[(i)] $\phi=[\alpha,\beta]$ is a commutator in $\diff_c(-1,+1)$ (resp. $\diff_+(\S^1)$)~;
\item[(ii)] The supports of $\psi_1$~, \dots, $\psi_r$ are two by two disjoint,
and their union brackets (resp. is nonempty and disjoint from) those of $\alpha$~, $\beta$~.
\end{enumerate}
Then $H$ is fillable without
 interior leaves.

\end{lem}
Except "without interior leaves", this is due to Thurston
\cite{thurston76}. Our construction differs from his,
in order to avoid interior compact leaves. We shall essentially need only the cases $r=1, 2$~.
\begin{proof}
{\def\middle{{\rm middle}}
\def\extremal{{\rm extremal}}

First the straight case: $Y=\D^1$~.
By condition (ii), there exist two points $-1<c_1<c_2<+1$~,
none of which is a common fix point to $\psi_1$~,\dots, $\psi_r$~, and between which
 the supports of
 $\alpha$ and $\beta$ lie.
So, there is a finite union $K\subset(c_1,c_2)$
 of compact subintervals
 such that the supports
of $\alpha$ and $\beta$ are contained in the interior of $K$~, and the supports of $\psi_1$~,\dots, $\psi_r$~, in $\D^1\setminus K$~.

One endows the solid cylinder $\D^2\times\D^1$ with
the Morse function $f$ obtained from the height function by creation of
 two singularities $s_1=(z_1,y_1)$~, $s_2=(z_2,y_2)$ of respective
 indices $1$~, $2$~,
in cancellation position. One arranges that $f(s_1)=c_1$ and $f(s_2)=c_2$~.
Stabilizing $f$ by
$\T^r$~, one gets the mapping:
$$F:H=\D^ 2\times\T^r\times\D^1\to\T^r\times\D^1:
(z,\theta,y)\mapsto(\theta,f(z,y))$$
One endows $\T^r\times\D^1$  with the suspension $\SS(\psi_1,\dots,\psi_r)$ of
 $\psi_1$~, \dots, $\psi_r$~.
Since $\psi_1,\dots,\psi_r$ are the identity on $K$~,
one can arrange that $\SS(\psi_1,\dots,\psi_r)$ is the height foliation in restriction to
 $\T^r\times K$~. Since $c_1, c_2$ are not common fix points to $\psi_1$~,\dots, $\psi_r$~, one can arrange that
$\SS(\psi_1,\dots,\psi_r)$ is transverse to
both tori $\T^r\times c_1$ and to $\T^r\times c_2$
 (the hypothesis that the supports of  $\psi_1$~, \dots, $\psi_r$
are two by two disjoint
makes this particuliarly immediate). 
Then, the mapping $F$ is transverse to the
 foliation $\SS(\psi_1,\dots,\psi_r)$~. Indeed, $F$ is a submersion, except
on both tori $z_1\times\T^r\times y_1$ and  $z_2\times\T^r\times y_2$~,
which $F$ respectively
maps diffeomorphically onto $\T^r\times c_1$ and $\T^r\times c_2$~,
both
transverse to $\SS(\psi_1,\dots,\psi_r)$~.
Consequently $F^*(\SS(\psi_1,\dots,\psi_r))$ is a (regular) foliation on $H$~.
(The description of this foliation, although not necessary for the proof,
is nevertheless of interest. Typically, for $r=1$~, one has a ``saddlization'':
 $z_1\times\S 1\times y_1$ for example, becomes the transverse intersection of
some 2-dimensional
Reeb component with some 3-dimensional Reeb component turning in the opposite
direction.)

 In $\D^2\times\D^1$~, the domain $f\mun(K)$
 is diffeomorphic to the product of $K$ with the compact
surface of genus one bounded by one circle. The fundamental group of this surface
being non-abelian free on two generators, the suspension of $\alpha$ and $\beta$ gives a
 foliation $\SS(\alpha,\beta)$ of this domain.

The foliation $\F$ filling $H$ is defined as two pieces: in $F\mun(\T^r\times K)$
one takes the suspension $\SS(\alpha,\beta)$ stabilized by the $r$-torus;
 and in the complement, $\F$ is  $F^*(\SS(\psi_1,\dots,\psi_r))$~.

Obviously, $\F$ fills the hole, i.e. coincides along $\partial H$
with the given germ.
Certainly $\F$ has
no leaf interior to $H$~. This is immediate in view of its trace
 on each $3$-dimensional
slice $\D^2\times\theta\times\D^1$~, for $\theta$ in $\T^r$~. 
Every leaf of this trace, being either
a level set of $f$ or a leaf of the suspension $\SS(\alpha,\beta)$~, meets the boundary
 $(\partial\D^2)\times\D^1$~.

\medbreak
The round case is much alike the straight one.
 By condition (ii)~, there exists  a finite union $K\subset\S^1$
 of compact subintervals
 such that the supports
of $\alpha$ and $\beta$ are contained in the interior of $K$~, and the supports of
 of $\psi_1,\dots,\psi_r$~, in $\S^1\setminus K$~.
By (ii) there are two points $c_1$~, $c_2$ on the circle,
both not common fixed points to $\psi_1,\dots,\psi_r$~, and that
lie in a same connected component of $\S^1\setminus K$~.
Arrange that $c_2<c_1$ in this interval endowed with the orientation induced from $\S^1$~.

Alike before, one endows the solid torus $\D^2\times\S^1$ with
the $\S^1$-valued Morse function $f$~, obtained from the
 height function, that is, the second projection, by creation of
 two singularities $s_1$~, $s_2$ of respective
 indices $1$~, $2$~,
in cancellation position. One arranges that $f(s_1)=c_1$
and $f(s_2)=c_2$~.
Thus the level set
$f\mun(y)$ has genus 1 (resp. 0) in case $y$ lies in the
connected component of $\S^1\setminus\{c_1,c_2\}$ containing (resp. not containing)
 $K$~.

Just as in the straight case, $F$ being the $r$-torus-stabilization of $f$~,
we fill $F\mun(\T^r\times K)$ by the
 $r$-torus-stabilized
suspension of $\alpha$ and $\beta$~; and the complement, by the suspension of
$\psi_1,\dots,\psi_r$ pulled back through $F$~.
}\end{proof}

\medbreak
Now let us come back to the construction of a minimal foliation
on $M$~. The cancellation of
the singularities
(paragraph \ref{cancel_sbs}) has left a minimal foliation $\xi_\holes$ with some standard, discal holes
of the form
 $H^n_i\cong\D^2\times\Sigma\times\D^1$~, where $1\le i\le n-2$ and $\Sigma=\S^{n-i-2}\times\S^{i-1}$~.
The monodromy is $\phi\neq\id$ over $\partial\D^2$~(\emph{compressible monodromy}),
and the identity over any loop in $\Sigma$~.
As already mentionned, $H^n_i$ is \emph{never} fillable.

However, $\xi_\holes$ being minimal, $M_\holes$ contains an embedded arc linking the
ceiling $\D^2\times\Sigma\times(+1)$ of $H^n_i$
to its floor $\D^2\times\Sigma\times(-1)$~, transversely to $\xi_\holes$~.
 One enlarges $H^n_i$ by a small tubular
neighborhood $W\cong\D^{n-1}\times\D^1$ of this arc, a
\emph{worm gallery}, obtaining a new hole $H_i^n\cup W$~. Its
germ of foliation along the portion $\S^{n-2}\times\D^1$ of its boundary, the one
that bounds $W$~, is of course the height foliation.
Actually, in the \emph{sub-extremal} cases $i=1$ or $n-2$~,
 the hole not being connected, we add two disjoint worm galleries:
one to each of the two connected components.

\begin{pro}\label{filling_pro} Let
 $n\ge 4$~, let $1\le i\le n-2$~, and let $\phi\in \diff_c(-1,+1)$~.
Then
the hole $H^n_i$ whose compressible monodromy is $\phi$~, enlarged by a worm gallery (or by two worm galleries
in case $i=1$ or $n-2$) is fillable
without interior leaf.
\end{pro}
Thurston fills the same holes with some kind of Reeb components inside \cite{thurston76}.
 We shall prove proposition \ref{filling_pro} in case $\phi$
 is a commutator $[\alpha_0,\beta_0]$ in $\diff_c(-1,+1)$~.

On the one hand, this case is enough to complete the proof of theorem A:
in the wirl previous to the cancellation of each pair of matched
 singularities (\ref{cancel_sbs}), take care to
choose for $\phi$ a commutator in $\diff_c(-1,+1)$~.
\def\LH{{\rm HOLED}}
Write $M_LH:=M$ minus the enlarged holes
and $\xi_\LH:=\xi\vert M_\LH$~, an obviously minimal foliation.
Applying
\ref{filling_pro}, extend $\xi_\LH$ inside each enlarged hole,
obtaining a foliation $\F$ on $M$~.
After the last words of proposition \ref{filling_pro},
 every leaf of $\F$
 meets $M_\LH$~, so its closure contains  $M_\LH$~. Every leaf being
locally dense, every leaf is dense in $M$~.

On the other hand,
our argument generalizes straightforwardly for a product of
 commutators,
and thus actually proves proposition \ref{filling_pro} in all cases,
 since  $\diff_c^\infty(\R)$
is simple (Thurston \cite{thurston74}, using previous results
by Epstein \cite{epstein70} and Herman
  \cite{herman}\cite{herman2}; see also \cite{epstein84}).
 But we don't need to use this difficult result.

\subsubsection{First
proof of proposition \ref{filling_pro} in dimension $4$}\label{dim_4_ssbs}
\begin{figure}
\includegraphics*[scale=0.6, angle=0]{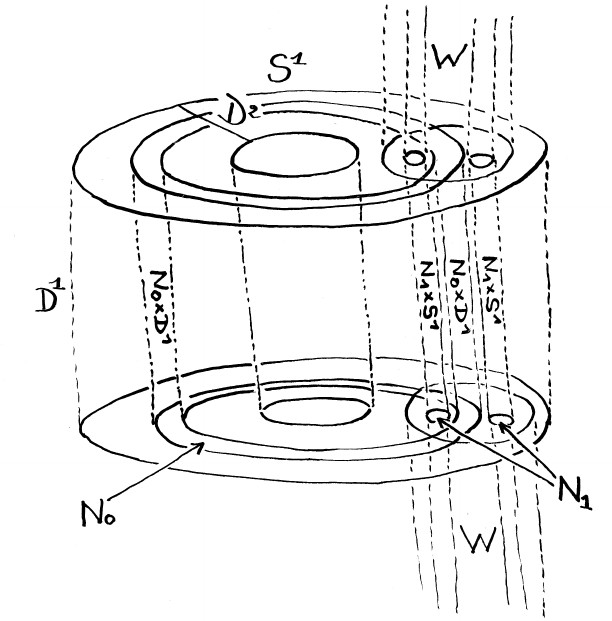}
\caption{The hole $H^4$ and its worm gallery.}
\label{hole_fig}
\end{figure}

\begin{figure}
\includegraphics*[scale=0.6, angle=0]{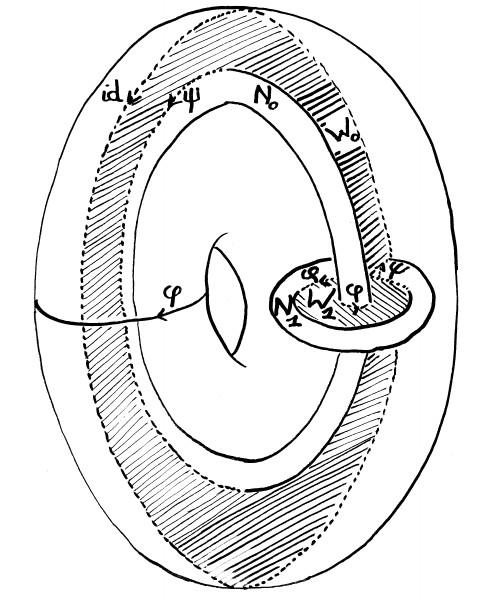}
\caption{Splitting of the base $B$ of the hole $H^4$~,
 and definition of the representation $\rho$~. The monodromy of a few loops is indicated;
$\psi$ is for $\phi_1$~.}
\label{base_fig}
\end{figure}

 In dimension $n=4$ there are two possible indices, $i=1$ or $2$~. The holes $H_1^4$ and $H_2^4$
are isomorphic by reversing the co-orientation of the foliation.
 Let $H^4$ be one of the two connected components of $H_1^4$~.
So, $H^4\cong\D^2\times\S^1\times\D^1$
is a straight discal hole whose base $B$ is
the solid torus $\D^2\times\S^1$~(figure
\ref{hole_fig}). Following Thurston,
we consider its core $V_0:=0\times\S^1$ and the meridian $V_1$ of $V_0$~;
and we split $B$ into three domains:
 two disjoint small compact tubular neighborhoods $N_0$~, $N_1$ of
$V_0$~, $V_1$~; and the complement  $C:=B\setminus Int(N_0\cup
N_1)$~.
In the base, $V_0$ bounds, modulo the boundary of the base,  the
annulus $W_0:=[0,1]\times\S^1$ (here $[0,1]$ is a radius
of the complex unit disk $\D^2$~); and $V_1$ bounds a 2-disk  $W_1$
that meets $V_0$ at its center (figure \ref{base_fig}).

Recall that $\phi$ is a
commutator $[\alpha_0,\beta_0]$ in
$\diff_c(-1,+1)$~, and
choose some $\alpha_1,\beta_1$ in $\diff_c(
-1,+1)$ whose supports are
disjoint from those of  $\alpha_0,\beta_0$~, and such that
the support of $\phi_1:=[\alpha_1,\beta_1]$
brackets those of $\alpha_0,\beta_0$ (see before lemma \ref{bracketting_inessential}).
In particular, $\phi_1$ commutes to $\phi$~.
One fills up $C\times\D^1$ with the suspension
 of the representation
$$\rho:\pi_1C\to \diff_c(-1,+1):\gamma\mapsto
\phi^{(W_0^*\gamma)}\phi_1^{(W_1^*\gamma)}$$
where $W_j^*\gamma$ denotes of course the algebraic intersection number between the loop $\gamma$ and the hypersurface $W_j$~. Obviously $\rho$ represents $\pi_1(\partial B)$ as needed:
 $\rho(\partial\D^2\times 1)=\phi$
and $\rho(1\times\S^1)=\id$~. It remains
to fill the two discal sub-holes $N_0\times\D^1$ and $N_1\times\D^1$ (union a worm gallery)~.

In restriction to
 $\partial N_0\cong
(\partial\D^2)\times\S^1$~, the representation $\rho$
associates $\phi$ to the factor $\partial\D^2$~, compressible
in $N_0$~;
and $\phi_1$ to the factor $\S^1$~, incompressible in $N_0$~.
 Since the support of $\phi_1$ brackets those
of $\alpha_0,\beta_0$~,
 the lemma \ref{bracketting_inessential}
applies, here $r=1$~, and fills the sub-hole $N_0\times\D^1$~.

On the contrary, in restriction to
 $\partial N_1\cong
(\partial\D^2)\times\S^1$~, the representation $\rho$
associates $\phi_1$ to the factor $\partial\D^2$~, compressible
in $N_1$~;
and $\phi$ to the factor $\S^1$~, incompressible in $N_1$~.
 So the lemma \ref{bracketting_inessential}
\emph{does not} apply to $N_1\times\D^1$~.

Now remember the worm gallery $W$~. The sub-torus $N_1$
is \emph{inessential} in $B$~, in the sense that $N_1$ is contained
in a 3-ball $\D^3$ embedded into the base. After an isotopy in $H^4$~,
one can assume that $\D^3\times\pm 1$ coincide with the entry and
the exit of the gallery. We extend the hole
$N_1\times\D^1$ through $W$~, obtaining an extended round hole $H_\round\cong N_1\times\S^1$~.
We keep the remaining of $W$ filled with the height foliation,
and so, the germ of foliation along $\partial H_\round$ is the suspension
of two commuting
diffeomorphisms $\bar\phi$~, $\bar\phi_1$ of the circle. Here
$\D^1$ is regarded as embedded into $\S^1$ and,
for every $\psi\in\diff_c(-1,+1)$~, one denotes $\bar\psi\in\diff_+(\S^1)$ its extension by the identity.
So, $\bar\phi_1=[\bar\alpha_1,\bar\beta_1]$~. The lemma \ref{bracketting_inessential}
applies to the round hole $N_1\times\S^1$~, and fills it, taking $r=1$~, taking
 $\bar\phi_1$ for $\phi$~, and taking $\bar\phi$ for $\psi_1$~.

Obviously, no interior leaf has been created in $H^4\cup W$~.
The first proof of proposition \ref{filling_pro}
in dimension 4 is complete.
\subsubsection{Proof of proposition \ref{filling_pro} for the subextremal indices}\label{subextremal_ssbs}
\medbreak
The preceding 4-dimensional construction is a pattern for the proof of proposition \ref{filling_pro} in all dimensions
$n\ge 4$~. The generalization is more direct in case the index of the hole is subextremal: $i=1$ or $n-2$~.
The hole $H^n_{n-2}$ is isomorphic to $H_1^n$~.
 Let $H^n$ be one of the
two connected components of $H^n_1$~, thus
a straight discal hole whose base is
 $B:=\D^2\times\S^{n-3}$~.

\begin{lem}\label{degree_lem}  For every $k\ge 0$~, there
is an embedding$$e_k:\S^{k}\times\S^1\to\Int(\D^ 1\times\S^{k+1})$$
such that the mapping $pr_2\circ e_k:\S^{k}\times\S^1\to\S^{k+1}$ is of degree one.
\end{lem}
For example, in case $k=0$~, the image $e_0(\S^{0}\times\S^1)$ consists of
two circles in the annulus, one essential, the other inessential. The general case is not much more difficult:

\begin{proof} One first embeds $\S^{k}\times\S^1$  into  $\S^{k+2}$ as the boundary
of a tubular neighborhood of a circle. Then $\D^1\times\S^{k+1}$ is also embedded
in   $\S^{k+2}$~, as the complement
 of two small $(k+2)$-balls, one interior to this neighborhood,
the other exterior to this neighborhood. Hence an embedding
$e_k$ from  $\S^{k}\times\S^1$ into $\D^1\times\S^{k+1}$~, whose compose
with the second projection is obviously of degree one.
\end{proof}

Here we depart from \cite{thurston76}, where one finds an analogous lemma, except that the $(k+1)$-torus
appears
where we prefer $\S^k\times\S^1$~. For the present task: filling the holes of subextremal indices, the difference
is anecdotic. In the following paragraph, where we shall fill the holes of intermediate indices, our choice will
simplify the construction: we shall have
to ``roll up the holes'' only twice.
But the critical difference
will be in the homotopy argument (section \ref{homotopy_sec}): it does not seem possible with tori. Simply,
$\S^{k-1}\times\S^1$ is in $\S^{k}\times\S^1$ the fix point set of an involution, while
$\T^{k}$ in $\T^{k+1}$ is not.
\medbreak
Now we begin to fill $H^n$~. Write $V:=\S^{n-4}\times\S^1$ and let
$j_0:V\to B$ be $e_{n-4}$ followed by the canonical
inclusion $\D^1\times\S^{n-3}\subset\D^2\times\S^{n-3}$~.
Put  $V_0:=j_0(V)$~. The mapping
 $$pr_2\circ j_0:V\to\S^{n-3}$$ being of degree one, $V_0$
together with $1\times\S^{n-3}$ bound
an orientable compact hypersurface $W_0$ in $B$ (recall that $1$ is a basepoint in $\partial\D^2$~).
 Of course, the bundle normal to $V_0$ in $B$ is trivial.
  One writes $X_0:=j_0(\S^{n-4}\times 1)$~; one writes
 $W_1$
 the total space of the (trivial) $\D^2$-bundle normal to $V_0$ above $X_0$~, so $W_1\cong\D^2\times\S^{n-4}$~;
 one writes
 $V_1:=\partial W_1$~; and one writes
 $N_0$~, $N_1$ two disjoint tubular neighborhoods of $V_0$~, $V_1$ in the interior of $B$~.

In $\D^1\times\S^{n-3}$~, every submanifold of
codimension two is inessential (see paragraph
 \ref{dim_4_ssbs}). In particular, $e_{n-4}(\S^{n-4}\times 1)$ is inessential
 --- this is also obvious from the
proof of lemma \ref{degree_lem}.
So, $X_0$ is inessential
in  $B$~.  So, $N_1$ is also.

One continues the filling of $H^n$ much as one did above (\ref{dim_4_ssbs})
 for $H^4$~: one assumes that $\phi=[\alpha_0,\beta_0]$ is a commutator;
 one chooses $\phi_1=[\alpha_1,\beta_1]$ just as one did in \ref{dim_4_ssbs}; 
one fills $C\times\D^1$~, where $C:=B\setminus Int(N_0\cup N_1)$~,
with the suspension of the representation $\rho$ defined as in \ref{dim_4_ssbs}.

There remains to fill two discal subholes $N_0\times\D^1$~, $N_1\times\D^1$~.

The first one $N_0\times\D^1$
 is just the $\S^{n-4}$-stabilization of a 4-dimensional
hole $\D^2\times\S^1\times\D^1$~. The lemma \ref{bracketting_inessential}, straight case, with $r=1$~, applies
to this 4-dimensional
hole and fills it with a foliation, whose  $\S^{n-4}$-stabilization
fills $N_0\times\D^1$~.

The second discal subhole $N_1\times\D^1$ is first extended, since $N_1$ is inessential in $B$~,
through the worm gallery. One
 obtains a round hole $N_1\times\S^1$~, which
 is just the $\S^{n-4}$-stabilization of a 4-dimensional
hole $\D^2\times\S^1\times\S^1$~.
 The lemma \ref{bracketting_inessential}, round case, with $r=1$~, applies
to this 4-dimensional
hole and fills it with a foliation, whose  $\S^{n-4}$-stabilization
fills $N_1\times\S^1$~.

 This completes the proof of proposition \ref{filling_pro} in every dimension $n\ge 4$~,
in case $i=1$ or $n-2$~. In dimension $4$~, this second proof of proposition \ref{filling_pro}
 will be the right one for the homotopy argument in the next section.

\subsubsection{Proof of proposition \ref{filling_pro} for the intermediate indices}\label{intermediate_ssbs}
In case $n\ge 5$ and $2\le i\le n-3$~, the base factor
 $\Sigma=\S^{n-i-2}\times\S^{i-1}$ is a
product of two spheres of positive dimensions. Unfortunately,
there exists, of course,
 no closed $(n-4)$-manifold $X$ that would
admit some mapping of nonnull degree$$X\times\S^1\to\S^{n-i-2}\times\S^{i-1}$$ that would be
null-homotopic in restriction to $X\times 1$~.
So we need one more iteration in the "rolling up" process: going somewhat to the meridian
of the meridian.
 
Consider the $(n-3)$-manifold
$$V:=\S^{n-i-3}\times\S^1\times\S^{i-2}\times\S^1$$
and the cartesian product $j_0:=e_{n-i-3}\times e_{i-2}$ of the two embeddings given by the lemma \ref{degree_lem}:
$$j_0: V\to\D^1\times\S^{n-i-2}\times\D^1\times\S^{i-1}\cong\D^2\times\Sigma=B$$
Put  $V_0:=j_0(V)$~. The mapping
 $pr_2\circ j_0:V\to\Sigma$ being of degree one, $V_0$
together with $1\times\Sigma$ bound
an orientable compact hypersurface $W_0$ in $B$~.
 Of course, the bundle normal to $V_0$ in $B$ is trivial.
  One defines
   $$X_0:=j_0(\S^{n-i-3}\times 1\times\S^{i-2}\times\S^1)$$
and one defines $W_1$ as
 the total space above $X_0$ of the (trivial) $\D^2$-bundle normal to $V_0$ in $B$~.
Let $V_1:=\partial W_1$ and
 choose a diffeomorphism $j_1:V\to V_1$ such that the normal projection of every point $j_1(x_1,\theta_1,
x_2,\theta_2)$ to $V_0$ is $j_0(x_1,1,
x_2,\theta_2)$~.
In the same way, one defines
$$X_1:=j_1(\S^{n-i-3}\times\S^1\times\S^{i-2}\times 1)$$
 and one defines $W_2$ as the
 total space over $X_1$ of the  $\D^2$-bundle normal to $V_1$ in $B$~. Let $V_2:=\partial W_2$~, a third embedding
of $V$ into the interior of $B$~; and let $N_0$~, $N_1$~, $N_2$ be two by two disjoint,
compact tubular neighborhoods of  $V_0$~, $V_1$~, $V_2$~, respectively, in the interior
of $B$~.

 This $N_2$ is inessential in $B$~.
Indeed, every submanifold of codimension 2 being inessential in the product of a sphere of positive dimension
 with an interval, the bisphere $j_0(\S^{n-i-3}\times 1\times\S^{i-2}\times 1)$ is inessential
in  $B$~.  But clearly, shrinking the fibre 2-disks over $X_0$ and $X_1$~, one brings $N_2$ arbitrarily close to this bisphere
through an isotopy in $B$~.

One continues the filling of $H_i^n$ ($2\le i\le n-3$) much as one did above
in the case $i=1, n-2$~.
In $\diff_c(-1,+1)$~, some diffeomorphisms $\alpha_j$~, $\beta_j$~,  $\phi_j$ ($0\le j\le 2$)
are chosen such that:
\begin{itemize}
\item $\phi_0=\phi$~;
\item Each $\phi_j=[\alpha_j,\beta_j]$~;
\item
For all $j\neq k$~, the supports of  $\alpha_j,\beta_j$ are disjoint from those of
 $\alpha_k,\beta_k$~;
 \item  For every $j\ge 1$~, the support of
 $\phi_j$
brackets those of $\alpha_{j-1}$ and $\beta_{j-1}$~.
\end{itemize}
One writes
 $$C:=B\setminus Int(N_0\cup N_1\cup N_2)$$
and fills the sub-hole $C\times\D^1$ by the suspension of the representation
$$\rho:
\gamma\mapsto\phi_0^{(W_0^*\gamma)}\phi_1^{(W_1^*\gamma)}\phi_2^{(W_2^*\gamma)}$$
There remains to fill
the sub-holes $N_0\times\D^1$~,  $N_1\times\D^1$ and  $N_2\times\D^1$~.

The first one is the stabilization by $X:=\S^{n-i-3}\times\S^{i-2}\times\S^1$ of a 4-dimensional hole $\D^2\times\S^1\times\D^1$~, whose monodromy
is $\phi_0$ over $\partial\D^2$ and $\phi_1$
over the incompressible $\S^1$ factor. The lemma \ref{bracketting_inessential} (with $r=1$) applies
to $\D^2\times\S^1\times\D^1$ and fills it with a foliation, whose $X$-stabilization fills 
 $N_0\times\D^1$~.
 
 The second hole $N_1\times\D^1$ is the $(\S^{n-i-3}\times\S^{i-2})$-stabilization of a 5-dimensional hole $\D^2\times\T^2\times\D^1$~, whose monodromy
is $\phi_1$ over $\partial\D^2$ and $\phi_0$~, $\phi_2$
over the two incompressible $\S^1$ factors. The lemma \ref{bracketting_inessential} (with $r=2$) applies
to $\D^2\times\T^2\times\D^1$ and fills it with a foliation, whose $(\S^{n-i-3}\times\S^{i-2})$-stabilization fills 
 $N_1\times\D^1$~.
 
  The third hole $N_2\times\D^1$ is the $X$-stabilization of a 4-dimensional hole $\D^2\times\S^1\times\D^1$~, whose monodromy
is $\phi_2$ over $\partial\D^2$ and $\phi_1$
over the incompressible $\S^1$ factor.
Since $N_2$ is inessential in $B$~,
one can extend this hole through the worm gallery, and one
 obtains a round hole $N_2\times\S^1$~, which
 is just the $X$-stabilization of a 4-dimensional
round hole $\D^2\times\S^1\times\S^1$~.
 The lemma \ref{bracketting_inessential} (with $r=1$) applies
to $\D^2\times\S^1\times\S^1$ and fills it with a foliation, whose $X$-stabilization fills 
 $N_2\times\S^1$~. This completes the proof of proposition \ref{filling_pro}
in the case $n\ge 5$ and $2\le i\le n-3$~,
and the construction of a minimal foliation on $M$~.

\section{Prescribing the homotopy class of the foliation as a $\G$-structure.}\label{homotopy_sec}

In this section, one more step is made towards theorem A:
given a closed connected manifold $M$ of dimension $n\ge 4$ whose Euler characteristic is null~, and
given on this manifold a smooth, co-oriented
 $\G$-structure $\xi$~,
 we build a homotopy from $\xi$ to some smooth, minimal foliation.

Essentially, the preceding section \ref{minimize_sec} already
did the job. Recall that it
 started with an arbitrary Morse foliation on $M$~, and then made it regular and minimal through some successive
transformations.
Now we shall more precisely start with a Morse foliation homotope to $\xi$~; and just have to
verify that each transformation is actually a homotopy
 of $\G$-structures.

They are in fact naturally so, because they verify some
 \emph{induction} property:
 the transformation is local, exists in every dimension $n$
from some minimal one, and is reflection-symmetric w.r.t.
some separating hypersurface, in restriction to which
it induces the same transformation in dimension $n-1$~.
So, conversely, one can extend the transformation from $M\times 1$
into $M\times[0,1]$ as one half of the same transformation
performed in dimension $n+1$~.
To verify that some transformation is inductive,
we shall most often content ourselves with pointing
the
induction hypersurface,
 the details then being straightforward; but a few points deserve some care.
\medbreak
As a basic example, our
two tools are inductive:

a) The creation
of a pair of Morse singularities $s$~, $s'$ of respective
 indices $i$~, $i+1$ is defined for $n\ge i+1$~, and inductive.
To fix ideas, regard it as a transformation of the height
foliation in $\R^n$~.
One can arrange that
the resulting function 
 is
 symmetric w.r.t. the hyperplane
 $0\times
\R^{n-1}$~, in restriction to which
one has the creation of a pair of singularities of \emph{the same} indices $i$ and $i+1$ (not to be confused,
 if $i\ge 1$~,
with the orthogonal hyperplane, in restriction to which
one has the creation of a pair of singularities of indices $i-1$ and $i$~.)

b) The cancellation of a pair of Morse singularities
 $s$~, $s'$ of respective
 indices $i$~, $i+1$~, verifying the hypotheses of lemma \ref{cancel_lem},
is also inductive for $0\le i\le n-1$~.
The induction
 hypersurface
 contains the unstable disk $D'$ and meets the stable disk $D$ transversely
into an $(n-i-1)$-disk.

\medbreak

Our aim is to make
some $\G$-structure $\eta$ on $\bar M:=M\times[0,1]$
such that $\eta\vert(M\times 0)=\xi$ and that
$\eta\vert(M\times 1)$ is regular and minimal. First consider the
$\G$-structure $\pr_1^*\xi$ on
$\bar M$~.
 Recall that,
as a $\G$-structure, $\pr_1^*\xi$ is nothing but a
 foliation $\X$ on some neighborhood
 of the
null section $\bar M$
 in the trivial bundle $\bar M\times\R$~, transverse
 to the $\R$-fibres. One endows $\bar M$ with
the $\G$-structure $\eta:=\sigma^*\X$ where
$$\sigma:\bar M\to\bar M\times\R$$is
a smooth section of this bundle,
 $C^0$-close to the null one,
 null over $M\times 0$~,
and such that $\sigma\vert(M\times 1)$
 has only quadratic tangencies with $\X$~.
 Since $\X$ is $C^2$~, that last condition is generic.
Then, $\eta\vert(M\times 1)$ is a Morse $\G$-structure
on $M\times 1$~. Call every point $s\in M\times 1$
that is singular for
$\eta\vert(M\times 1)$~, a
\emph{boundary singularity}.
By the index $i$ of $s$~,
we always mean its index in $M\times 1$~.
We turn $s$
into a \emph{half Morse singularity} of $\eta$~,
by a local, $C^0$-small perturbation of $\sigma$~, relative to $M\times 1$~. That is, after the perturbation,
 $\bar M$ admits local coordinates
$x_0=1-\pr_2$~, $x_1$~,\dots, $x_n$~, w.r.t. which
 $\eta$ is defined by the quadratic form
$$\pm x_0^2+x_1^2+\dots+x_{n-i}^2-x_{n-i+1}^2-\dots-x_n^2$$
We can choose the sign in front of $x_0^2$~. The convenient choices will be made later. Call $s$ {\it relatively stable} or
 {\it relatively unstable}
depending upon whether the sign is $+$
or $-$~.

Now we make $\eta\vert(M\times 1)$
 regular and minimal through the successive
transformations of the section \ref{minimize_sec}; and are going to verify that each transformation is inductive, and
thus extends from
$M\times 1$ to $\bar M$~, relatively to $M\times 0$~, as
one half of the same transformation
performed in dimension $n+1$~.

\medbreak
The replacement of a relatively stable
singularity $s$ of index $0$ by one of
index $2$ and a Reeb component,
as in paragraph \ref{minimize_sbs}, is defined in every dimension $n\ge 3$~,
and inductive. This follows
easily from a) and b).
The induction hypersurface
contains the singularities $s$~, $s'$~, $s"$
 and the arc $D'$~, and intersects the $n$-disk $D$
into a $(n-1)$-disk.

To apply this to some boundary singularity $s$ of index $0$~, first
we make it relatively stable; then, the transformation being inductive,
it extends to $\bar M$ as one half of the replacement of
a singularity of index $0$ in dimension $n+1$~.
 Finally, one has changed $s$ for a relatively stable boundary
singularity of index $2$ and a half $(n+1)$-dimensional Reeb component.

In the same way, every boundary
singularity of index $n$ is made relatively unstable, and then changed for a relatively unstable
boundary singularity of
index $n-2$ and a half $(n+1)$-dimensional Reeb component.

\medbreak
Also, the minimization process of paragraph \ref{minimize_sbs}
 is defined for
every $n\ge 3$ and inductive,
as follows from a).
The induction hypersurface
is $$\D^{n-2}\times\D^1\subset\D^{n-1}\times\D^1$$ 
The
four singularities of $f$ lie in
this hypersurface, as well as
their unstable invariant manifolds for some symmetric descending pseudo-gradient.
Then, the whirl is also inductive.

Thus, the minimization of $\eta\vert(M\times 1)$ extends into
$\bar M$ as a transformation of $\eta$~, relative to $M\times 0$~.
\medbreak
Then, as in \ref{cancel_sbs}, one creates new pairs of
boundary singularities for $\eta\vert(M\times 1)$~,
 until they match into pairs of successive indices.
After (a), this transformation extends to $\bar M$
as one half of the creation
of pairs of singularities in dimension $n+1$~.
 The new boundary singularities are relatively
stable.
\medbreak

The whirl that turns two small, disjoint, stable and unstable disks
of two singularities of indices $i$ and $i+1$~, into a cancellation pair, at the price of a hole  (paragraph \ref{cancel_sbs}),
is defined for $n\ge i+2$ and inductive. Here
the induction hypersurface is
$$\R^{n-i-2}\times\R^i\times\R\subset
\R^{n-i-1}\times\R^i\times\R$$

Now, consider a matched pair $s$~, $s'$ of boundary singularities
of $\eta$~, of indices $i$~, $i+1$~, with $1\le i\le n-2$~.
In case $i=n-2$~, reverse the co-orientation of $\eta$ to
shift to the case $i=1$~. So, $1\le i
\le n-3$~.
Make $s$ and $s'$ relatively stable. Then, the whirl that gives
them a cancellation pair, and their
 cancellation,
extend to $\bar M$ as one half of the
same transformation in dimension $n+1$~, and leaves
an unfoliated \emph{half hole}
$$\half H^{n+1}_i\cong
\D^2\times\half\S^{n-i-1}\times\S^{i-1}\times\D^1$$
whose intersection with $M\times 1$ is $H_i^n$~.
One writes $\half\D^k$ and $\half\S^k\cong\D^k$~, respectively, the compact half unit $k$-disk and the compact half unit $k$-sphere.
Thanks to the minimality of $\eta\vert(M\times 1)$~,
one  enlarges each connected component of $\half H^{n+1}_i$
with a \emph{half worm gallery} $\half W^{n+1}\cong\half\D^n\times\D^1$~,
whose intersection with $M\times 1$ is the worm gallery $W^n\cong
\D^{n-1}\times\D^1$~.
\medbreak
 There remains to fill the enlarged half holes,
that is, to verify that the filling of the enlarged holes,
performed in \ref{filling_sbs}, is inductive.
For
$i=1$ and for every $n\ge 5$~, in the enlarged hole
$$H_1^n\cup W^n=(\D^2\times\S^{n-3}\times\S^{0}
\times\D^1)\cup(\D^{n-1}\times\D^1)\cup(\D^{n-1}\times\D^1)$$
the induction hypersurface is the enlarged hole
$$H_1^{n-1}\cup W^{n-1}=(\D^2\times\S^{n-4}\times\S^{0}
\times\D^1)\cup(\D^{n-2}\times\D^1)\cup(\D^{n-2}\times\D^1)$$ 
For every
$i\ge 2$ and every $n\ge i+4$~, in the enlarged hole
$$H_i^n\cup W^n=(\D^2\times\S^{n-i-2}\times\S^{i-1}
\times\D^1)\cup(\D^{n-1}\times\D^1)$$
the induction hypersurface is the enlarged hole
$$H_i^{n-1}\cup W^{n-1}=(\D^2\times\S^{n-i-3}\times\S^{i-1}
\times\D^1)\cup(\D^{n-2}\times\D^1)$$
That is, the inclusion of the enlarged holes is induced by the inclusions $\S^{n-i-3}\subset\S^{n-i-2}$
and $\D^{n-2}\times\D^1\subset\D^{n-1}\times\D^1$~.
One has to verify that the filling of
$H_i^n\cup W^n$ induces
by restriction
the filling of
$H_i^{n-1}\cup W^{n-1}$~.
(This is the point that would fail if one had allowed
 $i=n-3$~, since
the number of galleries would not be the same for $H_i^{n-1}$
and for $H_i^n$~).

The inductivity of the hole filling is based on the following
straightforward inductive version of lemma \ref{degree_lem}.
Write $\S^\infty$ the inductive limit of
the spheres $\S^k$~.

\begin{lem}\label{inductive_degree_lem} There is a mapping
$$e:\S^{\infty}\times\S^1\to\D^1\times\S^{\infty}$$
such that
\begin{enumerate}
\item For each $k\ge 0$~, the restriction
 $e_k:=e\vert(\S^{k}\times\S^1)$ is an embedding
of $\S^{k}\times\S^1$
into $\Int(\D^ 1\times\S^{k+1})$~;
\item For each $k\ge 0$~,
the compose $\pr_2e_k$
is a mapping of degree one from $\S^{k}\times\S^1$ onto $\S^{k+1}$~.
\end{enumerate}
\end{lem}
\begin{proof}
One first embeds $\S^{\infty}\times\S^1$  into  $\S^{\infty}$ as the boundary
of a tubular neighborhood of $\S^1$~. Two points are chosen 
in $\S^{2}$,
one interior to this neighborhood, the other, exterior.
In $\S^\infty$~, the complement
 of two small balls centered at these two points is
homeomorphic to
 $\D^1\times\S^{\infty}$~,
and contains $\S^{\infty}\times\S^1$~, hence the embedding $e$~. 
Obviously, $\pr_2e\vert(\S^k\times\S^1)$
 is of degree one onto $\S^{k+1}$~.
\end{proof}

Now we apply the filling process of paragraphs \ref{subextremal_ssbs}
and \ref{intermediate_ssbs},
 using
the embeddings $e_k$ given by lemma
\ref{inductive_degree_lem} (For $n=4$ and $i=1$~,
we \emph{don't} apply paragraph \ref{dim_4_ssbs},
 whose construction is
not inductive, and of which I don't know if it builds
a homotopy of $\G$-structure.)

\emph{Subextremal index ---} Consider the case $i=1$~.
Since $e_{n-4}$ extends $e_{n-5}$~, in the base $B^{n-1}\cong\D^2\times\S^{n-3}$ of the connected component $H^n$~, one
easily arranges that $V_0$~, $W_0$~, $X_0$~, $W_1$~, $V_1$~, $N_0$~,
$N_1$~, $C$ intersect the base $B^{n-2}\cong\D^2\times\S^{n-4}$ of $H^{n-1}$ into the analogous subsets. Then, the representation
$\rho$ and the filling of $C$ are inductive.
The filling of $N_0$ by means of a stabilization of lemma
\ref{bracketting_inessential} is inductive, since the induction carries on
the stabilizing factors $\S^{n-5}\subset\S^{n-4}$~.
 In the same way,
$N_1$ being inessential in $B^{n-1}$ and $N_1\cap B^{n-2}$ being inessential in $B^{n-2}$~, one
chooses a $(n-1)$-ball $\D^{n-1}\subset B^{n-1}$ containing $N_1$
and intersecting $B^{n-2}$ into a $(n-2)$-ball $\D^{n-2}$~. One
 takes $\D^{n-1}$ (resp. $\D^{n-2}$)
 as the entrance and exit of the worm gallery attached to $H^n$
(resp. $H^{n-1}$). Then, the extension of $N_1\times\D^1$
through the gallery, as a round hole, is inductive. Its filling is inductive, for the same reasons as for $N_0$~. Thus, the filling
of the enlarged holes of index $1$ is inductive.

\emph{Intermediate index ---}
In case $2\le i\le n-4$~, completely similar arguments show
that the filling of the enlarged holes of intermediate index is also inductive.

\section{Prescribing the homotopy class of the foliation as a hyperplane field.}\label{hyperplane_sec}

In this section, the third and last step is made towards theorem A:
given a closed connected manifold $M$ of dimension $n\ge 4$~, and
given on this manifold a smooth, nonsingular
vector field $V$~,
 we build a smooth, minimal, co-oriented foliation whose descending gradient 
 is homotopic to $V$~.

Essentially, the section \ref{minimize_sec} already
did a large part
 of the job. Recall that at the end of
 paragraph \ref{minimize_sbs} one has
obtained a first minimal, co-oriented Morse foliation $\xi_\min$ on $M$~.
 Then (paragraph \ref{cancel_sbs}) one rather arbitrarily created pairs
of new singularities, so as to match the singularities
into pairs of successive indices; and one chose arbitrarily some arcs, transverse to $\xi$~, joining
the matched pairs. Later on, these arcs were used to cancel
the matched pairs, at the price of holes; and after filling
the holes, a minimal foliation was produced.
Now we shall show that if we are more careful in the
creation of the pairs of new singularities and in
the choice of the arcs, then the gradient of the produced foliation
 will
be naturally homotopic to $V$~.

One easily verifies that
this section and the previous one are compatible, allowing to prescribe in the same time
the homotopy class of the desired foliation as a $\G$-structure and the homotopy class
of its gradient, and thus establishing theorem A in the smooth, co-oriented case.

\medbreak
To avoid irrelevant technicallities, one fixes,
in a neighborhood of every singularity,
 a local system of
coordinates  $x_1,\dots,x_n$~, in which $\xi_\min$ is defined by
 $$
x_1^2+x_2^2+\dots-x_{n-i+1}^2-x_{n-i+2}^2\dots$$
Call them \emph{preferred}. In the rest of this section,
 by the words ``relative to the singularities''
 one means relative to \emph{some neighborhood} of the singularities.

 Using a partition of unity, one makes a
 \emph{pseudogradient} for $\xi_\min$~, that is,
 a vector field
 $\nabla\xi_\min$ negatively transverse to $\xi_\min$~, except at the
singularities of $\xi_\min$~, close to which one
has in the preferred coordinates
$$\nabla\xi_\min=-x_1\partial/\partial x_1-x_2\partial/\partial x_2-\dots$$$$+x_{n-i+1}\partial/\partial x_{n-i+1}
+x_{n-i+2}\partial/\partial x_{n-i+2}+\dots$$

After a homotopy of
 $V$~, one has, close to each singularity,
 $V=\partial/\partial x_1$ or $\partial/\partial x_{n-i+1}$~.
Call $V$ \emph{preferred}.

\medbreak

The argument lies on an (obvious) generalization of the
elementary
Thom-Pontryagin construction for the mappings $M\to\S^{n-1}$~.
The generalization is in the same time twisted and relative.
It concerns the nonsingular vector fields $X$
 on $M$ that coincide with $V$
on a neighborhood of $\sing(\xi_\min)$~.
To every such $X$~,
one associates the \emph{opposition} set $C_-(X,\nabla\xi_\min)$~:
the points of $M$ where $X$ and $\nabla\xi\min$ are nonpositively
colinear. If, over $M\setminus\sing(\xi_\min)$~, the
vector field $X$ is
transverse to $\nabla\xi_\min$ as sections of the projectivized tangent bundle, then $C_-(X,\nabla\xi_\min)$ is in $M$ a compact
submanifold of dimension $1$~, whose boundary is exactly $\sing
(\xi_\min)$~. Write$$C_-^*(X,\nabla\xi_\min)
:=C_-(X,\nabla\xi_\min)\setminus\sing(\xi_\min)$$
This curve is naturally
endowed with a $\tau\xi_\min$-framing $F(X,
\nabla\xi_\min)$~.
That is, $F(X,\nabla\xi_\min)$ is a vector bundle isomorphism,
at every point of $C_-^*(X,\nabla\xi_\min)$~,
from the $(n-1)$-space tangent to $\xi_\min$~, onto the
$(n-1)$-space
normal to $C_-(X,\nabla\xi_\min)$ in $M$~.

This way, the homotopy classes, rel. $\sing(\xi_\min)$~,
of the non\-sing\-ular vector fields that coincide with $V$
on a neighborhood of $\sing(\xi_\min)$~, are in one-to-one
correspondance with the
classes of
 $\tau\xi_\min$-framed 1-submanifolds bounded by $\sing(\xi_\min)$~.
The classes are taken w.r.t.
 $\tau\xi_\min$-framed
cobordism in $M\times[0,1]$ rel. $\sing(\xi_\min)\times[0,1]$~.
We don't need to think too much about the subtleties of this cobordism relation, because we shall anyway use only the following:
if$$C_-(X,\nabla\xi_\min)=C_-(V,\nabla\xi_\min)$$
and if 
$F(X,\nabla\xi_\min)$ is homotopic to $F(V,\nabla\xi_\min)$
rel. $\sing(\xi_\min)$~, then $X$ is homotopic to $V$
rel. $\sing(\xi_\min)$~.

\medbreak
After a small generic perturbation of $V$ relative to the singularities, $V$ and $\nabla\xi_\min$ are transverse as
sections of the projectivized tangent bundle.
One considers the projection $T$ of $V$
into $\tau\xi_\min$ parallelly to $\nabla\xi_\min$~. That is, $T$
is a smooth vector field on $M\setminus\sing(\xi_\min)$~, tangential to
$\xi$~, and one has a decomposition~:
 $$V=h\nabla\xi_\min+T$$ where $h$ is a function on $M\setminus\sing(\xi_\min)$~.
So, $T$ and $h$ have no common zero,
and
$C_-^*(V,\nabla\xi_\min)$
 is the set of
 zeroes of $T$ where $h<0$~.

One calls $x\in C_-^*(V,\nabla\xi_\min)$ \emph{nondegenerate}
if $T$~, regarded as a vector field in the leaf $L_x$
of $\xi$ through $x$~, admits at $x$ a nondegenerate zero.
Then, $C_-^*(V,\nabla\xi_\min)$ is transverse to $\xi$ at $x$~,
and
 the Poincar\'e-Hopf
index 
 of $T\vert L_x$ at $x$
 is $\pm1$~.

For example, close to each singularity $s$ of $\xi_\min$ of
Morse
index $i$~, because of the preferred form of $V$~,
every opposition point is nondegenerate, and
 the Poincar\'e-Hopf index
is $(-1)^i$ or $(-1)^{i+1}$ depending
upon whether $C_-(V,\nabla\xi_\min)$ is above
or below $s$~, w.r.t.
 the co-orientation
of $\xi_\min$~.

Also, note that 
at a nondegenerate $x$~, the
framing $F(V,
\nabla\xi_\min)$
is of course nothing but the 1-jet of $T$ in $L_x$~.

One calls a component of $C_-(V,\nabla\xi_\min)$
 \emph{nondegenerate} if every point of this component is nondegenerate, but the singularities of $\xi_\min$~.

By a \emph{bisingular arc,}
 we mean an arc embedded in $M$~,
 linking two singularities
of $\xi_\min$ whose indices are successive: $i$ and $i+1$~;
and transverse to $\xi_\min$ except at its extremities.

\begin{pro}\label{bisingular_pro} After the creation of pairs of singularities in $\xi_\min$
(not local extrema)
and after some homotopy of $V$ relative to the singularities,
$\xi_\min$ admits a pseudogradient $\nabla\xi_\min$
such that every connected component of $C_-(V,\nabla\xi_\min)$ is a nondegenerate bisingular
arc.
\end{pro}

To prove it, one first notices that
one has, in $M\setminus\sing(\xi_\min)$~, a neighborhood of $C^*_-(V,\nabla\xi_\min)$~,
such that any modification of $T$ in this neighborhood,
 relatively to
the singularities, into another vector field $T'$
also tangential to $\xi_\min$~,
 results into a homotopy of $V$~. Namely, if
the support of $T'-T$ is contained in
$\{h<0\}$~, then
for every $0\le t\le 1$~,
 $$V_t:=h\nabla\xi_\min+(1-t)T+tT'$$is nonsingular.
\medbreak
As a first application, we put $C_-(V,\nabla\xi_\min)$ in general
position w.r.t. $\xi_\min$~. Close to each point of $C^*_-(V,\nabla\xi_\min)$~,
 the field $T$~, viewed through any local chart trivializing the
foliation $\xi_\min$~,
 becomes a 1-parameter family of vector fields
on $\R^{n-1}$~.  By Shoshitaichvili's normal form
theorem \cite{sho72}
\cite{sho75},
after a small generic perturbation of $T$~, the bifurcations
are "saddle-nodes", or "folds". That is,
at every degenerate $x\in C^*_-(V,\nabla\xi_\min)$~,
one has the birth or the death
of a pair of singularities of $T$~. "Birth" is distinguished from
 "death"
by the co-orientation of $\xi_\min$~. In particular,
$C_-(V,\nabla\xi_\min)$
 has a quadratic tangency with $\xi_\min$ at $x$~,
and the Poincar\'e-Hopf
index vanishes at $x$ and changes sign. Call $x$ \emph{cubic.}

\medbreak
The toolbox to modify $C_-(V,\nabla\xi_\min)$ and its position w.r.t. $\xi_\min$
contains three tools.

Creation of a pair
of cubic points, a birth cubic point and a death cubic point.
This is easily achieved by changing
$T$ in a small neighborhood of any nondegenerate point of $C^*_-
(V,\nabla\xi_\min)$~.

Cancellation of a pair of cubic points. Let $x, x'$ be
on  $C_-^*(V,\nabla\xi_\min)$
 a death cubic point and a birth cubic point~;
and let $A$ be an arc in $M$~,
 from $x$ to $x'$~,
positively transverse to $\xi_\min$~, and otherwise disjoint from
$C^*_-(V,\nabla\xi_\min)$~. A first homotopy
of $V$~, supported in a small neighborhood of $A$~, makes
$h<0$ along $A$~, not changing $T$~. Then,
a second homotopy of $V$~, also
 supported in a small neighborhood of $A$~,
changing $T$ and not $h$~,
results, regarding $C_-$~, in
a surgery of index 1: a small subarc $\alpha$
of $C_-$ through $x$ and a small subarc $\alpha'$
of $C_-$ through $x'$
 are cut from $C_-$~, and two arcs parallel to $A$
are pasted, linking each endpoint of $\alpha$ to the
endpoint of $\alpha'$ with the same
Poincar\'e-Hopf index.

The third tool changes a cubic point for a pair of singularities of $\xi_\min$~.
This requires more care.

\begin{lem}\label{create_negative} Let $U\subset M$ be a nonempty
open subset, disjoint from $C_-(V,\nabla\xi_\min)$~,
 and let $0\le i\le n-2$~.

Then, there are
\begin{itemize}
\item A homotopy of $V$
with support in $U$~;
\item A homotopy of $\xi_\min$~, with support in $U$~,
that creates two singularities $s$~, $s'$
 of indices $i$~, $i+1$~;
\item A coherent
 modification of $\nabla
\xi_\min$ on a compact subset of $U$~;
\end{itemize}
whose effect on $C_-(V,\nabla\xi_\min)$ is to add a new arc component
 $C_-^{new}$
 linking $s$ to $s'$ in $U$~. All points on $C_-^{new\ *}$ are nondegenerate,
except one death cubic point.
\end{lem}

\begin{proof}{of lemma \ref{create_negative}.} 
 This will be verified on some local model where moreover
one is reduced to the 2-dimensional case (figure \ref{opposition_fig},
 left).
\begin{figure}
\includegraphics*[scale=0.5,angle=90]{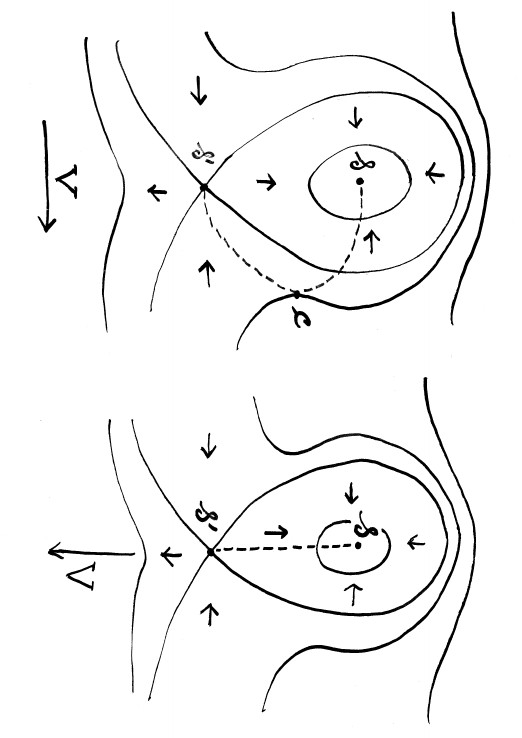}
\caption{Opposition curve
 of some constant vector field $V$ in the plane,
with some pseudogradient of a pair of singularities in cancellation position.
Two cases are represented. Dotted line: the opposition curve;
 $s, s'$~: singularities; $c$~: cubic point; arrows:
 descending pseudogradient.}
\label{opposition_fig}
\end{figure}
After passing to appropriate local coordinates, $U=\R^n$ endowed with coordinates $(x_1,\dots,x_n)$~, and in $U$~:

a) The original
$\G$-structure $\xi_\min$ is defined by a nonsingular
 function of the form~$$f:=x_n+Q(x_1,\dots,x_{n-1})$$where $Q=\lambda_1x_1^2+\dots+\lambda_{n-1}x_{n-1}^2$
 is a nondegenerate diagonal quadratic form
of rank $n-1$ and index $i$ (and $\lambda_1>0)$~;

b) The original pseudogradient is~:
$$\nabla\xi_\min=-(\partial f/\partial x_1,\dots,\partial f/\partial x_n)$$

Then, the vector field $V$ being in $U$ nowhere negatively colinear to $\nabla\xi_\min$~,
one can arrange moreover, after a homotopy with compact support, that
 in a neighborhood of $\D^{n-1}\times\D^1\subset\R^n$~,
 one has $V=
\partial/\partial x_1$~.

 Let $\xi$ be the $\G$-structure on $M$ that coincides with $\xi_\min$ outside $U$
and that admits in $U$ for a first integral the function~: 
$$g:=x_n-b(\Vert\hat x\Vert)b(x_n)+Q(x_1,\dots,x_{n-1})$$where $\Vert\hat x\Vert:=(x_1^2+\dots+x_{n-1}^2)^{1/2}$ and where $b:\R
\to[0,1]$ is a smooth
 bump function. Obviously, $\xi$ is homotopic to $\xi_\min$~.

 Define the pseudogradient $\nabla\xi$ as equal to $\nabla\xi_\min$ outside $U$~,
while inside $U$~:
$$\nabla\xi:=-(\partial g/\partial x_1,\dots,\partial g/\partial x_n)$$

Choose $b$ even, such that its support is exactly $[-1,+1]$~, such that $b=1$ in a neighborhood of $0$~, and such that $b''=0$
at exactly two points $\pm y_0$~. So, $b'=1$ at exactly two points $y_1$~, $y_2$~; and $-1<y_1<y_0<y_2<0$~.
Moreover choose the $\lambda_i$'s large enough with respect to the variations of the bump function,
namely:
$$\vert b'(y)/y\vert<2\vert \lambda_j\vert\ \ (y\in\R^*, 1\le j\le n-1)$$

 Then, the following is immediately verified.

The function $g$ has exactly two singularities $s:=(0,\dots,0,y_1)$~, $s':=(0,\dots,y_2)$~. Provided that $b$ has been chosen to coincide with convenient
polynomials of degree two in a neighborhood of $y_1$
and $y_2$~, both singularities are Morse, of indices $i$~, $i+1$~;
and indeed at $s$ and $s'$ the coordinates $
x_1$~,
\dots, $x_{n}$
 are preferred coordinates for $\xi$ (up to some permutation, additive constants,
and multiplicative constants).
(Some constant multiple of) $V$  is also preferred there w.r.t. these preferred coordinates.

The equations of the tangency curve $C(V,\nabla\xi)\cap U$ are$${\partial g\over\partial x_2}=\dots={\partial g\over\partial x_{n}}=0$$
This curve is contained in the 2-plane $x_2=\dots=x_{n-1}=0$~.
Its negative half $C^{new}_-:=C_-(V,\nabla\xi)\cap U$
is an arc linking the singularities.
 The equation of the degenerate points~:
$${\partial g\over\partial x_2}=\dots={\partial g\over\partial x_{n}}={\partial^2 g\over\partial x_n^2}=0$$
admits a unique solution $c:=(t_0,0,\dots,y_0)$ on $C^{new}_-$~.
For a generic choice of $b$ in a neighborhood of $t_0$ and $y_0$~, this unique degenerate point $c$ is a death cubic point.
\end{proof}
Here is our third tool. One is given on $C_-(V,\nabla\xi_\min)$ a cubic point $c$~, 
say, a birth point;
and an index $1\le i\le n-2$~. One creates in $M$~,
 close to $c$
and under it (w.r.t. the co-orientation of $\xi_\min$),
a pair of singularities of $\xi_\min$ bounding a new arc component
of $C_-(V,\nabla\xi_\min)$~, on which there is a unique, death, cubic
point $c'$ (lemma \ref{create_negative}).
 A short arc, transverse to $\xi_\min$~,
is chosen from $c'$ to $c$~. Then, $c$ and $c'$ are cancelled
(tool 2). The result on $C_-(V,\nabla\xi_\min)$ is that a small
subarc $\alpha$ through $c$ has been cut from $C_-(V,\nabla\xi_\min)$~, and that
the extremity of $\alpha$ whose
 Poincar\'e-Hopf index was $(-1)^i$ (resp. $(-1)^{i+1}$)
 has become a singularity of $\xi_\min$
of index $i$ (resp. $i+1$~).

Symmetrically, given a death cubic point $c$~, a small
subarc $\alpha$ through $c$ can be cut from $C_-(V,\nabla\xi_\min)$~, and the
extremity of $\alpha$ whose
 Poincar\'e-Hopf index was $(-1)^i$ (resp. $(-1)^{i+1}$)
 becomes a singularity of $\xi_\min$
of index $i+1$ (resp. $i$~).

\begin{proof}{ of proposition \ref{bisingular_pro}}.
After having  used tool 1 to create
 a pair
of cubic points on every circle component of $C_-(V,\nabla\xi_\min)$~, change every cubic point on $C_-(V,\nabla\xi_\min)$
into a pair of singularities (tool 3). Then, every component $A$
is an arc free of cubic points, in particular transverse to $\xi$~. Orienting $A$ by the co-orientation of $\xi$~, let $s$ (resp. $s'$)
be the lower (resp. upper) endpoint,
 $i$ (resp. $i'$) its Morse index,
 and $\delta(A):=i'-i$~. It remains to force $\delta(A)=1$~.

We have seen that
 the Poincar\'e-Hopf degree, which is constant
along $A$~, equals $(-1)^i$ close to $s$~,
and $(-1)^{i'-1}$ close to $s'$~. So, $\delta(A)$ is odd.

If $\delta(A)\ge 3$~, one creates a pair of cubic points on $A$ (tool 1), and one changes
each of them
into a pair
of singularities of indices
$i+1$ and $i+2$ (tool 3). So, $A$ has been cut into
$3$ arcs without cubic points. For the first and the second, $\delta=1$~;
while the third $A'$ has
$\delta(A')=\delta(A)-2$~.

If $\delta(A)\le -1$~, one creates a pair of cubic points on $A$~.
One changes the first one from $s$ into a pair
of singularities of indices
$i+1$ and $i$~; and the second, into a pair
of singularities of indices
$i-1$ and $i-2$~. So, $A$ has been cut into
$3$ arcs without cubic points. For the first and the second,
 $\delta=1$~;
while the third $A'$ has
$\delta(A')=\delta(A)+2$~.

By successive applications of these two cases, one makes $\delta=1$ for
all components of $C_-(V,\nabla\xi_\min)$~.
\end{proof}

Then, one proceeds much as in paragraph \ref{cancel_sbs}: for each component $A$
of $C_-(V,\nabla\xi_\min)$~,
two small disks $D$~, $D'$~, stable and unstable,
 are chosen at its extremities
$s$~, $s'$~, such that $A\cap D$ is a radius of $D$
and $A\cap D'$ is a radius of $D'$~. One chooses, as
in paragraph \ref{cancel_sbs}, an embedding $e:\R^n\to M$
whose image is a neighborhood
$U$ of $A\setminus(A\cap(D\cup D'))$~,
in the same position w.r.t. $\xi_\min$ ,
$D$ and $D'$~, as in paragraph \ref{cancel_sbs}, but with the extra property
that $e\mun(A)$ is the $x_n$-axis. The whirl of $\xi_\min$
in each $U$ (that is,
applied simultaneously close to all connected components
of $C_-(V,\nabla\xi_\min)$)
 gives on $M$ a $\G$-structure with holes $\xi_\holes$ that admits
at $s, s'$ a cancellation pair $\Phi(D), D'$~. One obviously can
perform
the whirls so that $\xi_\holes=\xi_\min$ on a small neighborhood of
$C_-(V,\nabla\xi_\min)$~. One fixes
a pseudogradient $\nabla\xi_\holes$ for $\xi_\holes$~, equal
to $-\partial/\partial x_n$ on each $U$~. The term
"pseudogradient" for a \emph{holed} $\G$-structure means a vector
field on $M$ which is a pseudogradient on $M_\holes$~,
and the height gradient inside the holes. The virtue
of $\nabla\xi_\holes$ is that it is a common pseudogradient for $\xi_\min$ and $\xi_\holes$~.

Ideally, a more precise version of lemma \ref{cancel_lem}
would then produce a foliation with a pseudogradient
 whose opposition curve with
$\nabla\xi_\holes$ is $A$~.
The following is enough for our purpose, and avoids delicate
computations.

\begin{lem}\label{cancel_lem_precise}
 Let
\begin{itemize}
\item $M$~,
 $\xi$~, $s$~, $s'$~, $D$~, $D'$ be as in lemma \ref{cancel_lem};
\item $N$
be a small enough neighborhood of $D\cup D'$ in $M$~;
\item $A$ be the bisingular arc between $s$ and $s'$~,
 union of a radius of $D$
with a radius of $D'$~;
\item some preferred coordinates
be fixed for $\xi$ at $s,s'$~;
\item $\nabla\xi$ be a pseudogradient for $\xi$ in $N$~.
\end{itemize}
Then, $N$
 admits a foliation $\NN$ and a nonsingular vector field $X_N$~,
such that
\begin{itemize}
\item One has $\NN=\xi$
and $X_N=\nabla\xi$~, but on some compact subset in $N$~;
\item $X_N$ is homotopic to some pseudogradient $\nabla\NN$ of $\NN$
rel. the complement of some compact subset of $N$~;
\item $X_N$ is
preferred at $s$ and $s'$~;
\item
$C_-(X_N,\nabla\xi)$ coincides with $A$ and is nondegenerate; 
\item No leaf of $\NN$
is relatively compact in $N$~.
\end{itemize}
\end{lem}
\begin{proof} One follows the proof of lemma
\ref{cancel_lem}, with some more
precision and some differences.
If $N$ is small enough, $\xi$ admits in $N$
 a first integral $f$~, null
 in restriction to the bouquet of spheres $\partial D\cup\partial D'$~.
This bouquet has a compact neighborhood
$E\subset N$ diffeomorphic to $B\times[-\epsilon,
+\epsilon]$~, where $B$ is
 $\S^{n-i-1}\times\S^i$ minus
an open $(n-1)$-disk, and such that
 $f\vert E$ is the height function.

Then, changing $f$ on some smaller neighborhood of $D\cup D'$~,
makes the critical value at $s$ lift up to $-\epsilon/2$~,
and
 the critical value at $s'$ descend down to $\epsilon/2$~.
One gets on $N$ a Morse function $f'$~, equal to $f$ outside
some compact subset, and whose singularities
are still $s$ and $s'$~. It is well-known, and easy to verify,
 that one can choose $f'$ so
that $\nabla\xi$ is still a pseudogradient for $f'$~,
and so that $f'$ defines $\xi$
on a neighborhood of $D\cup D'$~.

 The bouquet $\partial D\cup\partial D'$ has
a compact neighborhood $E'\subset N$ diffeomorphic to
$\D^{n-1}\times[-\epsilon,+\epsilon]$~; such that
$f'\vert E'$ is the height function on a neighborhood of
 $\partial E'$~. The singularities $s$~, $s'$
 are in cancellation position in $E'$~, in
Morse's sense: one easily makes a pseudogradient for $f'\vert E'$
which is the height gradient close to $\partial E'$~, and
tangential and radial on $D$ and on $D'$~. In particular, the
stable invariant manifold of $s$ and the unstable invariant
manifold of $s'$ meet transversely, and their intersection is
$A$~.

On the other hand, there is, in  $\D^{n-1}\times\D^1$~,
 a model function with a cancellable pair of singularities,
for which we can control its opposition curve with the height gradient (figure \ref{opposition_fig}, right). Endow $\D^{n-1}\times\D^1$
with the coordinates $x_1$~,\dots, $x_n$ (so, $x_n$
is the height function);  with the function $g$ from lemma
\ref{create_negative};
 and with the pseudogradient
$$\nabla g:=-(\partial g/\partial x_1,\dots,
\partial g/\partial x_n)$$
Write $\sigma,\sigma'$ the two singularities of $g$~.
 Recall that
 $
x_1$~,
\dots, $x_{n}$
 are preferred coordinates for $g$ at $\sigma,\sigma'$ (up to some permutation, additive constants,
and multiplicative constants).
(Some constant multiple of) $-\partial/\partial x_n$  is also preferred there, w.r.t. these preferred coordinates.
One verifies easily, by the same computations
made to prove lemma \ref{create_negative}, that $C_-(-\partial/\partial x_n,
\nabla g)$ is the straight segment $[\sigma,\sigma']=0\times[y_1,y_2]$~, and nondegenerate.
The singularities $\sigma$~, $\sigma'$
 are in cancellation position in $\D^{n-1}\times\D^1$~, in
Morse's sense.

By elementary Morse theory, $g$ is conjugate to $f'$~,
through some diffeomorphism $\zeta:\D^{n-1}\times\D^1\to E'$~.
One can choose $\zeta$ so that $\zeta[\sigma,\sigma']=A$~;
and so that the preferred
coordinates for $\xi$ at $s, s'$
as the $\zeta$-image of those for $g$ at $\sigma,\sigma'$~.

One
defines $\NN$ by $x_n\circ\zeta\mun$
 on $E'$~, 
and by $f'$ on $\R^n\setminus E'$~.
One defines $\nabla f'$ on $N$ as a second
pseudogradient for $f'$~, equal to
 $\nabla\xi$ outside some compact subset,
and to $\zeta_*\nabla g$
on $E'$~.
One defines a pseudogradient $\nabla\NN$ for $\NN$
as $-\zeta_*\partial/\partial x_n$ on $E'$~,
and $\nabla f'$ on $N\setminus E'$~.

So, $C_-(\nabla\NN,\nabla f')=A$~, and it
is nondegenerate. Consider in $N\setminus\{s,s'\}$
 the decomposition $\nabla\NN=
h\nabla f'+T$~, where $T$ is tangential to $f'$~.
Recall that $\nabla\xi$ is also a pseudogradient for $f'$~.
One has on $N$ a convex homotopy,
with compact support, of nonsingular vector fields
$$X_t:=h((1-t)\nabla f'+t\nabla\xi)+T$$
carrying $X_0=\nabla\NN$ to $X_N:=X_1=h\nabla\xi+T$~,
whose opposition curve with $\nabla\xi$ is $A$~.
Decompose $X_N$ as $h'\nabla\xi+T'$~, where $T'$
is tangential to $\xi$~. Since we arranged that $f'$
defines $\xi$ close to $A$~, one has $T=T'$ there. So,
the opposition $C_-(X_N,\nabla\xi)$ is nondegenerate.
\end{proof}

This lemma, applied simultaneously to all the components $A$
of $C_-(V,\xi_\min)$~, gives on $M$~:
\begin{itemize}
\item
A minimal foliation with holes $\F_\holes$~, equal to $\NN$ on each $N$
and to $\xi_\holes$ on the complement;
\item
A pseudogradient  $\nabla\F_\holes$~, equal to $\nabla\NN$
 on each $N$
and to $\nabla\xi_\holes$ on the complement;
\item
A nonsingular vector field $X$~,
 equal to $X_N$ on each $N$
and to $\nabla\xi_\holes$ on the complement. So, $X$ is
homotopic to $\nabla\F_\holes$~, preferred at the singularities,
and $C_-(X,\nabla\xi_\holes)=C_-(V,\nabla\xi_\min)$
is nondegenerate.
\end{itemize}

Let again $A$ be any
component of $C_-(X,\nabla\xi_\holes)=C_-(V,\nabla\xi_\min)$~.
 Recall that, the ambiant dimension $n$ being
at least $4$~,
there are two homotopy
classes of framings along $A$~, rel. $\partial A$~.
Of course, there is no reason that $F(X,\nabla\xi_\holes)$ and $F(V,\nabla\xi_\min)$ belong to the same class.
If they don't, one makes a twist around $A$~. Namely,
 recall that one has chosen an embedding $e:
\R^n\to M$~. One changes $e$ for
its compose
with a twist of $\R^n$ around its $x_n$-axis, e.g. 
$$
(x_1,\dots,x_n)\mapsto
(\alpha^{\psi(x_n)}(x_1,\dots,x_{n-1}),x_n)$$
$\psi$ being a smooth function on $\R$
such that
$\psi(x)=0$ for $x\le-1/4$ and $\psi(x)=2\pi$ for $x\ge 1/4$~, and $(\alpha^t)$ being some $2\pi$-periodic, 1-parameter subgroup of
$\SO(n-1)$~. (Regarding the compatibility of the
present section with the section \ref{homotopy_sec}, this is
the only point that deserves the notice that it \emph{is}
compatible: such a twist is clearly inductive in the sense
of section \ref{homotopy_sec}.)
This change in the choice of $e$ amounts to twist $\xi_\holes$~,
$\nabla\xi_\holes$~, $\F$~, $\nabla\F$~, $X$~, and $
F(X,\nabla\xi_\holes)$~,
once around $A$  (the holes are twisted too).
After applying this twist around each $A$ where it is necessary,
$F(X,\nabla\xi_\holes)$ and $F(V,\nabla\xi_\min)$ are
homotopic w.r.t. $\sing(\xi_\min)$~. That is, $X$ has
the same Thom-Pontryagin invariant, w.r.t. $\nabla\xi_\holes$~,
as $V$ w.r.t. $\nabla\xi_\min$~. There remains a little homotopy
argument before we can apply Thom-Pontryagin.

Recall the decomposition$$V=h\nabla\xi_\min+T$$
$T$ being tangential to $\xi_\min$~.
One has a convex homotopy of vector fields
$$V_t:=h((1-t)\nabla\xi_\min+t\nabla\xi_\holes)+T$$
carrying $V_0=V$ to $V_1=h\nabla\xi_\holes+T$~. Every $V_t$
is nonsingular, since $\nabla\xi_\min$ and $\nabla\xi_\holes$ are
two pseudogradients for $\xi_\min$~.
Finally, decompose $V_1$ under the form
$$V_1=h_1\nabla\xi_\holes+T_1$$
$T_1$ being tangential to $\xi_\holes$~. Recall that,
close to $C_-(V,\nabla\xi_\min)$~, one has
$\xi_\min=\xi_\holes$~.
So, $T=T_1$ there. Finally,
$$C_-(V_1,\nabla\xi_\holes)=C_-(V,\nabla\xi_\min)=
C_-(X,\nabla\xi_\holes)$$
$$F(V_1,\nabla\xi_\holes)=F(V,\nabla\xi_\min)\sim
F(X,\nabla\xi_\holes)$$

By Thom-Pontryagin,  $X$
is homotopic to $V_1$~. That is, $\nabla\F_\holes$
is homotopic to $V$~.
\medbreak
There remains to fill the holes, as in paragraph \ref{filling_sbs}, so that
the produced foliation $\F$ has $\nabla\F\sim\nabla\F_\holes$~.

Obviously, $\F_\holes$ admits a pseudogradient which is the height
gradient not only in the holes $H_i^n$~, but also in their worm
galleries. Then, in the parts of the enlarged holes that one fills
by suspensions, one defines of course $\nabla\F$ as the height
gradient. So, one is reduced to arrange that, in lemma
\ref{bracketting_inessential},
the filling foliation $\F$ admit a pseudogradient $\nabla\F$
equal to the height gradient $-\partial/\partial y$
along $\partial H$~,
and homotopic to $-\partial/\partial y$ rel. $\partial H$~.
To fix ideas, consider the straight
 case. In paragraphs \ref{dim_4_ssbs}, \ref{subextremal_ssbs}
 and \ref{intermediate_ssbs} we can
 choose $\phi_1$ such that moreover,  $\phi_1(y)\ge y$
for every $y\in\D^1$~. So, in the hypotheses of lemma
\ref{bracketting_inessential},
 one has moreover $\psi_1(y)\ge y$~.
 Then, in $H$~,
 obviously $\F$ admits a pseudogradient $\nabla\F$~, equal to $-\partial/\partial y$ along $\partial H$~, such that
$d\theta_1\nabla\F\ge 0$ at every point of $H$~. Here, $\theta_1$
is, of course,
 the projection of $\T^r$ onto its first $\S^1$ factor.
Then, $\nabla\F$ and $-\partial/\partial y$~, being both nowhere
opposed to $\partial/\partial\theta_1$~, are homotopic rel.
$\partial H$~.

The round case is alike, the condition $\psi_1(y)\ge y$ referring to the connected component of $\S^1\setminus K$ containing $y$~,
oriented by the orientation of $\S^1$~.
The proof of Theorem A is complete.

\section{Generalizations}\label{generalizations_sec}
 Thurston's $h$-principle is indeed
 more general than the form given in the introduction, in three ways. Our theorem A
admits the same generalizations.
The compact manifold $M$ may have a smooth boundary.
 To simplify,
we consider only the construction of foliations \emph{transverse}
 to
this boundary. One works in 
 any differentiability class
$C^r$~, $1\le r\le\infty$~. The $\G$-structures are not
necessarily co-oriented.

In full generality, on the manifold $M$~,
a {\it Haefliger structure} of codimension one --- more briefly a \emph{$\G$-structure}  \cite{haefliger} ---
 can be defined as a pair $(\nu\xi,\xi)$~:
\begin{itemize}
\item A rank-one real
 vector bundle $\nu\xi$ over $M$~,
the \emph{normal} bundle;
\item In the total space
of $\nu\xi$~, a germ, along the zero section $Z(M)$~,
of foliation $\xi$
 of codimension one, transverse to every fibre.
\end{itemize}
A ($\nu\xi$)-\emph{twisted} tangent vector $V$ at point $x\in M$
is a linear morphism $\nu\xi_x\to\tau M_x$~.
If $V$ is transverse to $\xi$~,
the sign of this transversality is well-defined.
If $\xi$ is regular, it admits a \emph{gradient,}
that is, a $\nu\xi$-twisted vector field negatively
transverse to $\xi$~. 
 A homotopy between two $\G$-structures on $M$~, sharing the same normal bundle, is a $\G$-structure
on $M\times[0,1]$~.

\begin{TAP} Let
\begin{itemize}
\item $M$ be a compact manifold of dimension at least $4$~;
\item $\xi$
be a $\G$-structure of class $C^r$ ($1\le r\le\infty$)
 on $M$~, transverse to $\partial M$~;
\item $V$ be a $\nu\xi$-twisted
vector field on $M$~, tangential to $\partial M$~,
and a gradient of
$\xi\vert\partial M$~.
\end{itemize}
 Then, there is a minimal $C^r$ foliation
on $M$~, $C^r$-homotopic
 to $\xi$ rel. $\partial M$ (as a $\G$-structure),
 and whose gradient is homotopic to $V$ rel. $\partial M$
(as a nonsingular twisted vector field).
\end{TAP}
The generalization of the $h$-principle
 to manifolds with boundary is important for
the homotopic theory of foliations:
up to easy considerations on homotopy classes of hyperplane
fields, it reduces
 the
classification of foliations on closed manifolds, up
to concordance, to the classification of $\G$-structures,
 up to homotopy, for which Haefliger has built a classifying space.

As an example of application of theorem A' for manifolds with
boundary, any two linear foliations of codimension one on the
$3$-torus are concordant through some smooth
 foliation of $\T^3\times[0,1]$
\emph{without interior leaf}.

\medbreak
We indicate the few changes to make in 
sections \ref{minimize_sec} through
\ref{hyperplane_sec}
 to prove theorem A'.

To pass to the relative version, just change
$M$ for $\Int(M)$  wherever necessary:
one never has to modify $\xi$~, nor $V$~, on the boundary.
\medbreak
A lower differentiability class $C^r$  ($1\le r<\infty$)
makes few problems, essentially because one works in a small
neighborhood of a subcomplex of dimension 1 in $M$~:
 one can make
$\xi$ smooth (that is, $C^\infty$) there by some $C^r$-small homotopy before working.
One works with Morse $\G$-structures of class $C^r$~,
by which we
mean that the $\G$-structure is $C^r$ on $M$~, smooth on a neighborhood
of its singularities, and that its singularities are Morse.
All pseudogradients are understood smooth, as well as $V$~.
Two points call for an argument, especially in the $C^1$ case.

In section
\ref{homotopy_sec}, the foliation $\X$ is only
$C^r$~; and we need a $C^0$-small section
$\sigma$ of $\nu\xi$ over $M\times 1$~,
such that $\sigma^*\X\vert(M\times 1)$ is a Morse
 $\G$-structure of class $C^r$~.
One makes first a $C^0$-small, $C^1$ section $\sigma_0$
whose tangency points with $\X$ are isolated (\cite{lm}, proposition
2.1). Then, $\X$
is easily smoothen in a neighborhood
of these tangency points,
through some ambiant, $C^r$-small, isotopy in the total
space of $\nu\xi$. Finally,
one takes for $\sigma$ a generic smooth section, $C^1$-
close enough to $\sigma_0$~.

In section \ref{hyperplane_sec}, in the beginning of the proof of proposition 4.1,
we have to put $C_-(V,\nabla\xi_\min)$ in general
position w.r.t. the Morse $\G$-structure $\xi_\min$ of class  $C^r$~.  
 In a 1-parameter generic
family of vector fields of class $C^r$~, the 
degenerate singularities are isolated. So, after a generic,
$C^r$-small perturbation of $T$~, the opposition curve
$C_-(V,\nabla\xi_\min)$ is tangent to $\xi$ at only
finitely many points. Then, $\xi_\min$
is easily smoothen in a neighborhood
of this curve,
through some ambiant, $C^1$-small, $C^r$ isotopy in $M$~, rel. $\sing(\xi_\min)$~. The pseudogradient $\nabla\xi_\min$
is not changed. Then we can apply Shoshitaichvili's theorem.

 The rest
of the proof of theorem A' takes place in this neighborhood where
$\xi_\min$ is smooth, except for the transverse arcs that are
guidelines for the worm galleries (paragraph \ref{filling_sbs},
just before proposition \ref{filling_pro}). 
Before digging each gallery, one smoothens $\xi_\holes$ in
a small neighborhood of the arc, through some $C^1$-small,
ambiant isotopy of class  $C^r$~,
 supported in a small neighborhood of the arc.
\medbreak
In the same way, passing from the co-oriented case
 to the twisted case
makes no problem, essentially because one works in a
domain where $\nu\xi$ is orientable. More precisely,
at the beginning of the proof of proposition \ref{bisingular_pro},
$C_-(V,\nabla\xi_\min)$ becomes a disjoint union of arcs, in a
neighborhood of which $\nu\xi_\min=\nu\xi$ is thus orientable.
 The rest of the proof
of theorem A' takes place in this neighborhood, except, once
again,  for the 
worm galleries. They are small neighborhoods
of transverse arcs linking the ceiling of each hole
to its floor. So,
 the co-orientation of $\nu\xi_\min$
in the holes obviously extends into the galleries.

\end{document}